

\documentstyle{amsppt}

\define\M{{\Cal M}}
\define\R{{\Cal R}}
\define\E{{\Cal E}}
\define\La{{\Cal L}}
\define\Aa{{\goth A}}
\define\bb{{\goth B}}
\define\cc{{\goth C}}
\define\dd{{\goth D}}
\define\ee{{\goth E}}
\define\nem{\overline {NE}}
\magnification1200
\topmatter

\title
On the Picard number of Fano 3-folds\\ with terminal
singularities
\endtitle

\author
Viacheslav V. Nikulin
\endauthor

\address
Steklov Mathematical Institute,
ul. Vavilova 42, Moscow 117966, GSP-1, Russia.
\endaddress

\email
slava\@nikulin.mian.su
\endemail

\abstract
We prove:
Let $X$ be a Fano 3-fold with terminal $\Bbb Q$-factorial
singularities and $X$ does not have a small extremal ray
and a face of Kodaira dimension $1$ or $2$ for Mori
polyhedron $\nem (X)$. Then Picard number $\rho (X) \le 7$.
\endabstract

\endtopmatter

\document

\head
INTRODUCTION.
\endhead

Here we continue investigations started in \cite{N6}, \cite{N7}.

Algebraic varieties we consider are defined over the field
$\Bbb C$ of complex numbers.

In this paper, we get a final result on estimating of the
Picard number
$\rho = \dim N_1(X)$ of a Fano 3-fold $X$ with terminal
$\Bbb Q$-factorial singularities if $X$ does not have small extremal
rays and its Mori polyhedron does not have
faces with Kodaira dimension 1 or 2. One can consider this class as
a generalization of the class of Fano 3-folds with Picard number 1. There are
many non-singular Fano 3-folds satisfying this condition and with
Picard number 2 (see \cite{Mo--Mu} and also
\cite{Ma}). We also think that studying of  the
Picard number of this class may be important for studying of Fano
3-folds with Picard number 1 too (see Corollary 2 below).

Let $X$ be a Fano 3-fold with $\Bbb Q$-factorial terminal
singularities. Let $R$ be an extremal ray of the
Mori polyhedron $\nem (X)$
of $X$. We say that $R$ has the {\it type (I)} (respectively {\it (II)})
if curves
 of $R$ fill an irreducible divisor $D(R)$ of $X$ and the contraction
 of the ray $R$ contracts the divisor $D(R)$ in a point
(respectively on a curve).  An extremal ray $R$ is called {\it small} if
curves of this ray fill a curve on $X$.

A pair
$\{ R_1, R_2 \}$ of extremal rays has the type
$\bb_2$ if extremal rays $R_1, R_2$ are different, both
have the type (II), and have the
same divisor $D(R_1)=D(R_2)$.

We recall that a face $\gamma $ of Mori polyhedron
$\nem(X)$ defines a contraction $f_\gamma : X \to X^\prime$
(see \cite{Ka1} and \cite{Sh}) such that $f(C)$ is a point
for an irreducible curve $C$ iff $C$ belongs to $\gamma$.
The $\dim X^\prime$ is called the Kodaira dimension of the $\gamma $.
A set $\E$ of extremal rays is called extremal if it is contained in
a face of Mori polyhedron.

\proclaim{Basic Theorem} Let $X$ be a Fano 3-fold with terminal
$\Bbb Q$-factorial singularities. Assume that $X$ does not have
a small extremal
ray, and Mori polyhedron $\nem (X)$ does not have a face of Kodaira
dimension 1 or 2.

Then we have the following statements for the $X$:

(1) The $X$ does not have a pair of extremal rays of the type $\bb_2$
and Mori polyhedron $\nem (X)$ is simplicial;

(2) The $X$ does not have more than one extremal ray
of the type (I).

(3) If $\E$ is an extremal set of $k$ extremal rays of $X$, then
the $\E$ has one of the types:
$\Aa_1\amalg (k-1)\cc_1$, $\dd_2\amalg (k-2)\cc_1$,
$\cc_2\amalg (k-2)\cc_1$,
$k\cc_1$ (we use notation of Theorem 2.3.3).

(4) We have the inequality for the Picard number of the $X$:
$$
\rho (X) = \dim N_1(X) \le 7.
$$

\endproclaim

\demo{Proof} See Theorem 2.5.8.
\enddemo

It follows from (4):

\proclaim{Corollary 1} Let $X$ be a Fano 3-fold with terminal
$\Bbb Q$-factorial singularities and $\rho (X)>7$.

Then $X$ has either a small extremal ray or a face of Kodaira
dimension 1 or 2 for Mori polyhedron.
\endproclaim

We mention that non-singular Fano 3-folds do not have a small extremal
ray (by Mori \cite{Mo1}),
and the maximal Picard number for them is equal to 10 according to
their classification by Mori and Mukai
\cite{Mo-Mu}. Thus, all these statements already
work for non-singular Fano 3-folds.

 From the statement (2) of the Theorem, we also get the following application
of Basic Theorem to geometry of Fano 3-folds.

Let us consider a Fano 3-fold $X$ and blow-ups $X_p$ in different
non-singular points $\{x_1,...,x_p\}$ of $X$. We say that this is
a Fano blow-up if $X_p$ is Fano.  We have the following very simple

\proclaim{Proposition}
Let $X$ be a Fano 3-fold with terminal $\Bbb Q$-factorial
singularities and without small extremal rays.
Let $X_p$ be a Fano blow up of $X$.

Then for any small extremal ray $S$ on $X_p$, the
$S$ has a non-empty intersection with one of exceptional
divisors $E_1,...,E_p$ of this blow up and does not belong to any of
them. The divisors $E_1,...,E_p$ define $p$ extremal rays of the type
$(I)$ on $X_p$.
\endproclaim

\demo{Proof} See Proposition 2.5.14.
\enddemo

It is known that a contraction of a face of
Kodaira dimension 1 or 2 of $\nem (X)$ of a Fano 3-fold $X$ has
a general fiber which is a rational surface or curve respectively,
because this contraction has relatively negative canonical class.
See \cite{Ka1}, \cite{Sh}.
It is known that a small extremal ray is rational \cite{Mo2}.

Then, using Basic Theorem and Proposition, we can divide
Fano 3-folds of Basic Theorem on the following 3 classes:

\proclaim{Corollary 2}
Let $X$ be a Fano 3-fold with terminal $\Bbb Q$-factorial
singularities and without small
extremal rays, and without faces of Kodaira dimension 1 or 2 for
Mori polyhedron. Let $\epsilon$ be the number of extremal rays
of the type (I) on $X$ (by Basic Theorem, the $\epsilon \le 1$).

Then there exists $p$, $1\le p \le 2-\epsilon$, such that
$X$ belongs to one of classes (A), (B) or (C) below:

(A) There exists a Fano blow-up $X_p$ of $X$ with a face of Kodaira
dimension 1 or 2. Thus, birationally, $X$ is a fibration on
rational surfaces over a curve or rational curves over a surface.

(B) There exist Fano blow-ups $X_p$ of $X$ for
general $p$ points on $X$ such that
for all these blow-ups the
$X_p$ has a small extremal ray $S$. Then images of curves of $S$ on
$X$ give a system of rational curves on $X$ which cover a Zariski
open subset of $X$.

(C) There do not exist Fano blow-ups $X_p$ of $X$ for general
$p$ points.

We remark that for Fano 3-folds with Picard number $1$ the
 $\epsilon =0$. Thus, $1\le p \le 2$.
\endproclaim

Using statements (2), (3) and (4)
of Basic Theorem, one can formulate similar
results for Fano blow ups in curves.

\smallpagebreak

To prove Basic Theorem, we classify appropriate
so called extremal sets and
E-sets of extremal rays of the type
(I) or (II).  We use so called diagram method to deduce from
this classification the statement (4) of the Basic Theorem.

A set $\Cal E$ of extremal rays is called {\it extremal} if it is
 contained
in a face of Mori polyhedron. The $\Cal E$ has
{\it Kodaira dimension 3}
if a contraction of this face gives a morphism on a 3-fold.
For Fano 3-folds with
$\Bbb Q$-factorial terminal singularities,
we give a description of extremal sets $\E$
of Kodaira dimension 3
which contain extremal rays of the types (I) or (II) only.

A set $\Cal L$ of extremal rays is called {\it $E$-set}
if $\Cal L$ is not
 extremal, but any proper subset of $\Cal L$ is extremal. Thus, the
$\La$ is minimal non-extremal.
For Fano 3-folds with $\Bbb Q$-factorial terminal singularities,
we give a description of
$E$-sets $\La$ such that $\La$ contains
extremal rays of the types (I) or (II) only, and any proper subset of
$\Cal L$ is extremal of Kodaira
dimension 3.

I am grateful to Profs. Sh. Ishii, M. Reid and J. Wi{\'s}niewski
for useful discussions. I am grateful to referee for useful comments.

I am grateful to Professor Igor R. Shafarevich for his constant interest to
and support of these my studies.

This paper was prepared in Steklov Mathematical Institute, Moscow;
Max-Planck Institut f\"ur Matematik, Bonn, 1990; Kyoto University,
1992--1993; Mathematical Sciences Research Institute, Berkeley, 1993.
I thank  these Institutes for hospitality.

Preliminary variant of this paper was published as a preprint \cite{N8}
(see also \cite{N9} about connected results).

\newpage

\head
CHAPTER 1. DIAGRAM METHOD.
\endhead

Here we give a simplest variant of the
diagram method for multi-dimensional algebraic varieties.
We shall use this method in the next
chapter. This part also contains some corrections and generalizations
to the corresponding part of our papers \cite{N6} and \cite{N7}.

Let $X$ be a projective algebraic variety  with $\Bbb Q$-factorial
singularities over an
algebraically closed field. Let $\dim X \ge  2$. Let
$N_1(X)$ be the $\Bbb R$-linear space generated by all
algebraic  curves on  $X$ by the numerical equivalence,
 and let $N^1(X)$ be the  $\Bbb  R$-linear space generated
by all Cartier (or Weil) divisors  on X by the  numerical
equivalence. Linear spaces $N_1(X)$ and $N^1(X)$ are dual to one
another by the intersection pairing.
Let $NE(X)$ be a convex cone in $N_1(X)$ generated by
all effective curves on $X$.
Let $\overline {NE}(X)$ be the closure of the
cone $NE(X)$  in $N_1(X)$. It  is called
{\it Mori cone ({\rm or} polyhedron)} of $X$.
A non-zero element  $x \in N^1(X)$ is called  $nef$ if
$x\cdot \overline {NE}(X) \ge 0$. Let $NEF(X)$  be  the
set of all nef elements of $X$ and the zero.  It  is the
convex cone in $N^1(X)$  dual to Mori cone $\overline
{NE}(X)$.  A ray  $R\subset \overline {NE}(X)$ with origin 0
is called {\it extremal}  if  from $C_1 \in \overline
{NE}(X)$, $C_2\in \overline {NE}(X)$  and $C_1+C_2\in R$
it  follows that $C_1  \in R$ and $C_2    \in R$.

We consider a condition (i) for
a set $\Cal R$ of extremal rays on X.

(i) \it If $R \in \Cal R$, then all curves $C \in R$ fill out
an irreducible  divisor $D(R)$  on $X$.

\rm
In this case, we can correspond to $\Cal R$ (and subsets of $\Cal
R$) an
oriented graph $G({\Cal R})$ in the following way: Two different
rays $R_1$ and
$R_2$ are joined by an arrow $R_1R_2$ with the beginning in $R_1$
 and the
end  in $R_2$ if $R_1\cdot D(R_2)>0$. Here and in what follows,
for an extremal
ray  $R$ and a divisor $D$ we  write $R\cdot D>0$  if
$r\cdot D>0$
for $r\in R$ and $r\not= 0$. (The same for the symbols $\le$, $\ge$
 and $<$.)

A set $\Cal E$ of extremal rays is called {\it extremal}
if it is contained in  a
face of ${\overline {NE}(X)}$. Equivalenty, there exists a nef
element $H \in N^1(X)$ such that ${\Cal E}\cdot H=0$. Evidently,
a subset of an extremal set is extremal too.

We consider the following condition (ii) for extremal sets
$\E$ of extremal rays.

(ii) \it An extremal set
$\E = \{ R_1,...,R_n \}$ satisfies the condition
(i), and for any real numbers $m_1\ge 0,....,m_n\ge 0$ which are not
all equal to $0$,
there exists a ray $R_j\in \Cal E$ such that
$R_j\cdot (m_1D(R_1)+m_2D(R_2)+...+m_nD(R_n))<0$. In particular,
the effective divisor $m_1D(R_1)+m_2D(R_2)+...+m_nD(R_n)$ is not
$nef$.

\rm
A set $\La $ of extremal rays is called {\it $E$-set} (extremal in a
different sense) if the $\Cal L$ is
not extremal but every proper subset of $\Cal L$ is extremal.
Thus, $\La$ is a minimal non-extremal set of extremal rays.
Evidently, an $E$-set $\La$ contains at least two elements.

We consider the following condition (iii) for $E$-sets
$\La $.

(iii) \it Any proper subset of an $E$-set
$\La=\{ Q_1,...,Q_m\}$ satisfies the
condition (ii), and
there exists a non-zero effective nef divisor
$D({\Cal L})=a_1D(Q_1)+a_2D(Q_2)+...+a_mD(Q_m)$.

\rm
The following statement is very important.

\proclaim{Lemma  1.1}
An $E$-set $\La$ satisfying the condition (iii) is  connected in the
following sense: For any decomposition
${\Cal L}={\Cal L}_1\coprod {\Cal L}_2$, where
${\Cal L}_1$ and ${\Cal L}_2$ are non-empty, there  exists an
arrow  $Q_1Q_2$ such that $Q_1 \in {\Cal L}_1$ and $Q_2 \in
{\Cal L}_2$.

If ${\Cal L}$ and ${\Cal M}$ are two different
$E$-sets satisfying the condition (iii), then there exists an arrow
 $LM$ where
$L\in {\Cal L}$ and $M\in {\Cal M}$.
\endproclaim

\demo{Proof} Let $\La=\{ Q_1,...,Q_m\}$. By (iii), there
exists a nef divisor $D({\Cal L})=a_1D(Q_1)+a_2D(Q_2)+...+a_mD(Q_m)$.
If one of the coefficients $a_1,...,a_m$ is equal to zero, we get
 a contradiction with the conditions (ii) and (iii). It follows
that all the coefficients $a_1,...,a_m$ are positive. Let
${\Cal L}={\Cal L}_1\coprod {\Cal L}_2$
where
$\La_1=\lbrace Q_1,...,Q_k\rbrace$ and $\La_2=\lbrace Q_{k+1},...,Q_m
\rbrace $. The divisors
$D_1=a_1D(Q_1)+...+a_kD(Q_k)$ and $D_2=a_{k+1}D(Q_{k+1})+...+a_mD(Q_m)$
are non-zero.
By (ii), there exists a ray $Q_i, 1\le i\le k$, such that
$Q_i\cdot D_1<0$. On the other hand,
$Q_i\cdot D({\Cal L})=Q_1\cdot (D_1+D_2)\ge 0$. It follows, that
there exists $j$, $k+1\le j\le m$, such that $Q_i\cdot D(Q_j)>0$.
 It means that $Q_iQ_j$ is an arrow.

Let us prove the second statement. By the condition (iii), for every
 ray $R\in {\Cal L}$, we have the inequality
$R\cdot D({\Cal M})\ge 0$. If $R\cdot D({\Cal M})=0$
for any $R\in {\Cal L}$, then the set $\Cal L$ is
extremal, and we get the contradiction. It follows that there exists
 a ray $R \in {\Cal L}$ such that $R\cdot D({\Cal M})>0$.
It follows the statement.
\enddemo

Let $NEF(X)=\overline {NE}(X)^\ast \subset N^1(X)$ be the cone of nef
elements of X and ${\Cal M}(X)=NEF(X)/{\Bbb R}^+$ its projectivization.
We use usual relations of orthogonality between subsets of
$\M(X)$ and $\nem(X)$. So, for
$U\subset \M(X)$ and $V\subset \nem(X)$ we write
$U\perp V$ if
$x\cdot y =0$ for any ${\Bbb R}^+x\in U$ and any $y\in V$. Thus,
for $U \subset \M(X)$, $V\subset \nem(X)$ we denote
$$
U^\perp = \{ y\in \nem(X) \mid U\perp y\},\ \
V^\perp = \{ x\in \M(X)\mid x\perp V\}.
$$

A subset $\gamma \subset \M(X)$ is called
a {\it face} of $\M(X)$ if there
exists a non-zero element $r \in \nem(X)$ such that
$\gamma =r^\perp$.

A convex set is called a {\it closed polyhedron} if it is a
convex hull of a finite set of points. A convex closed polyhedron is called
{\it simplicial} if all its faces are simplexes. A  convex
closed polyhedron is called {\it simple} (equivalently,
it has simplicial angles) if it is dual to a simlicial one.
In other words, any its face of codimension
 $k$ is contained exactly in $k$ faces of $\gamma$ of
the highest dimension.

We need some relative notions to notions above.

We say that ${\Cal M}(X)$ is a {\it closed polyhedron in its face}
 $\gamma \subset {\Cal M}(X)$ if $\gamma $ is a closed polyhedron and
$\M(X)$ is a closed polyhedron
in a neighbourhood $T$ of $\gamma $.
Thus, there should exist a closed polyhedron
$\M^\prime$ such that $\M^\prime \cap T = \M(X) \cap T$.

We will use the following notation.
Let ${\Cal R}(X)$ be the set of all extremal rays of $X$.
For a face $\gamma \subset \M(X)$,
$$
{\Cal R}(\gamma )=\lbrace R\in {\Cal R}(X){ }\mid {  } \exists
{\Bbb R}^+H\in \gamma :R\cdot H=0\rbrace
$$
and
$$
\R(\gamma ^\perp )=\{ R \in \R(X) \mid \gamma \perp R \}.
$$
Let us assume that $\M(X)$ is a closed polyhedron in its face
$\gamma$.  Then sets
$\R(\gamma_1)$ and $\R(\gamma_1^\perp)$ are finite for any face
$\gamma_1 \subset \gamma $. Evidently, the face $\gamma $ is simple if
$$
\sharp \R(\gamma_1^\perp ) - \sharp \R(\gamma ^\perp )=
\text{codim}_\gamma \gamma_1
\tag1
$$
for any face $\gamma_1$ of $\gamma $. Then we say that the
polyhedron $\M(X)$ is {\it simple in its face $\gamma$}. Evidently,
this condition is equivalent to the condition:
$$
\dim~[\E] -\dim~[\R(\gamma^\perp)]=\sharp \E - \sharp \R(\gamma^\perp)
\tag2
$$
for any extremal set $\E$ such that
$R(\gamma^\perp)\subset \E$. Here
$[\cdot]$ denotes a linear hull.  (In \cite{N6}, we required a more strong
condition for a polyhedron $\M(X)$ to be simple in its face $\gamma$:
$\sharp  \R(\gamma_1^\perp )=\dim \M(X)-\dim \gamma_1$ for any face
$\gamma_1$ of $\gamma$.)

Let $A, B$ are two vertices of an oriented graph $G$. The
{\it distance} $\rho (A,B)$ in $G$ is a length (the number of links)
of a shortest oriented path of
the graph $G$ with the beginning in $A$ and
the end in $B$. The distance is $+\infty $ if this path does not
exist.  The {\it diameter} diam $G$ of an oriented graph $G$ is
the maximum distance between
ordered pairs of its
vertices. By the Lemma 1.1, the diameter of an $E$-set is a finite
number if this set satisfies the condition (iii).

The Theorem 1.2 below is an analog for algebraic varieties of
arbitrary dimension of the Lemma 3.4 from \cite{N2} and the Lemma 1.4.1
from \cite{N5}, which were devoted to surfaces.

\proclaim{Theorem 1.2}
Let $X$ be a projective algebraic variety
with $\Bbb Q$-factorial singularities and $\dim$ $X\ge2$. Let us suppose that
${\Cal M}(X)$
is closed and simple in its face $\gamma$.
Assume that the set ${\Cal R}(\gamma)$
satisfies the condition (i) above. Assume
that there are some
constants $d, C_1, C_2$ such that the conditions (a) and (b)
below hold:

(a)
For any $E$-set
$\La \subset {\Cal R}(\gamma )$ such that $\La$ contains at
least two elements
which don't belong to $\R (\gamma^\perp )$ and for any proper subset
$\La^\prime \subset \La$ the set $R(\gamma^\perp)\cup \La^\prime$ is
extremal, the condition (iii) is valid and
$$
\text{diam}~G({\Cal L})\le d.
$$

(b) For any extremal subset $\E$ such that
$\R(\gamma^\perp )\subset \E \subset \R(\gamma )$,
we have: the $\E$ satisfies the condition (ii)  and for the distance in
the oriented graph $G(\E)$
$$
\sharp \lbrace (R_1, R_2)\in (\E-\R(\gamma^\perp)) \times
(\E-\R(\gamma^\perp ))
\mid 1 \le \rho (R_1,R_2)\le d\rbrace
\le C_1\sharp ({\Cal E}-\R(\gamma^\perp ));
$$
and
$$
\split
\sharp \lbrace (R_1, R_2)\in (\E-\R(\gamma^\perp ))\times
(\E-\R(\gamma ^\perp))
\mid d+1\le \rho (R_1,R_2) &\le 2d+1\rbrace \\
                           &\le C_2\sharp ({\Cal E}-\R(\gamma^\perp )).
\endsplit
$$

Then $\dim ~\gamma <(16/3)C_1+4C_2+6$.

\endproclaim

\demo{Proof}
 We use the following Lemma 1.3 which was proved in \cite{N1}. The Lemma
 was
used in \cite{N1} to get a bound ($\le 9$) on the dimension of a hyperbolic
(Lobachevsky) space admitting an
action of an arithmetic reflection group with a field
of definition of the degree $>N$. Here $N$ is some absolute constant.

\proclaim{Lemma 1.3}
Let ${\Cal M}$ be a convex closed simple
polyhedron of a dimension $n$, and $A_n^{i,k}$ the average number of
$i$-dimensional faces of $k$-dimensional faces of $\Cal M$.

Then for $n\ge 2k-1$

$$
A_n^{i,k}< {{n-i\choose n-k}\cdot
({[n/2]\choose i}+{n-[n/2]\choose i})\over
{[n/2]\choose k}+{n-[n/2]\choose k}}.
$$
In particular, if $n\ge 3$
$$A_n^{0,2} <
\cases
4(n-1)\over n-2  &\text{if $n$ is even,}\\
4n\over n-1 &\text{if $n$ is odd.}
\endcases
$$
\endproclaim

\demo{Proof} See \cite{N1}. We mention that the right side of the
inequality of the Lemma 1.3 decreases and tends to the number
$2^{k-i}{k\choose i}$ of $i$-dimensional faces of $k$-dimensional
cube if $n$ increases.
\enddemo

 From the estimate of $A_n^{0,2}$ of the Lemma, it follows the following
analog of Vinberg's Lemma from \cite{V}. Vinberg's Lemma was used by
him to obtain an estimate ($\dim <30$) for the dimension
of a hyperbolic space admitting an action of
a discrete reflection group with a bounded fundamental polyhedron.

By definition, an {\it angle} of a polyhedron $T$ is an angle of a
2-dimensional face of $T$. Thus, the angle is defined by
a vertex $A$ of $T$, a plane containing $A$ and a 2-dimensional face
$\gamma_2$ of $T$, and two rays with the beginning at $A$ which
contain two corresponding sides of the $\gamma_2$. To
define an {\it oriented angle} of $T$, one should in addition
put in order two rays of the angle.

\proclaim{Lemma 1.4} Let $\Cal M$ be a convex simple polyhedron of
a dimension $n$. Let $C$ and $D$ are some numbers.
Suppose that oriented
angles ($2$-dimensional, plane) of $\Cal M$ are supplied with
weights and the following
conditions (1) and (2) hold:

(1) The sum of weights of all oriented angles at any vertex of
$M$ is not greater than $Cn+D$.

(2) The sum of weights of all oriented angles of any $2$-dimensional
face of $\Cal M$ is at least $5-k$ where $k$ is the number of
vertices of the $2$-dimensional face.

Then
$$
n<8C+5+
\cases
1+8D/n &\text{if $n$ is even,}\\
(8C+8D)/(n-1) &\text{if $n$ is odd}
\endcases .
$$
In particular, for $C\ge 0$ and $D=0$, we have
$$
n<8C+6.
$$
\endproclaim

\demo{Proof} We correspond to a non-oriented plane angle of $\Cal M$ a
weight which is equal to the sum of weights of
two corresponding oriented angles.
 Evidently, the conditions of the Lemma hold for the weights
of non-oriented angles too if we forget about the word "oriented".
Then we obtain Vinberg's lemma from \cite{V} which we
formulate a little bit more precisely here.
Since the proof is simple, we give the proof here.

Let $\Sigma$ be the sum of weights of all
(non-oriented) angles of the polyhedron
$\Cal M$. Let $\alpha _0$ be the number of vertices of
$\Cal M$ and $\alpha _2$ the number of $2$-dimensional faces
of $\Cal M$. Since $\Cal M$ is simple,
$$
\alpha _0{n(n-1)\over 2}=\alpha _2 A_n^{0,2}.
$$
 From this equality and conditions of the Lemma, we get inequalities
$$
(Cn+D)\alpha _0\ge \Sigma \ge \sum \alpha _{2,k}(5-k)=
5\alpha _2-\alpha _2A_n^{0,2}=
$$
$$
=\alpha _2(5-A_n^{0,2})=
\alpha_0(n(n-1)/2)(5/A_n^{0,2}-1).
$$
Here $\alpha _{2,k}$ is the number of 2-dimensional faces  with
$k$ vertices of $\Cal M$.  Thus, from this inequality
and Lemma 1.3, we get
$$
Cn+D\ge (n(n-1)/2)(5/A_n^{0,2}-1)>
\cases
n(n-6)/8 &\text{if $n$ is even,}\\
(n-1)(n-5)/8 &\text{if $n$ is odd.}
\endcases
$$
 From this calculations, Lemma 1.4 follows.
\enddemo

The proof of Theorem 1.2. (Compare with \cite{V}.)
Let $\angle$ be an
oriented angle of $\gamma $.
Let $\R(\angle)\subset \R(\gamma )$ be
the set of all extremal rays of $\M(X)$ which are orthogonal
to the vertex of $\angle$.
Since $\M(X)$ is simple in $\gamma$, the set
$\R(\angle)$ is a disjoint union
$$
\R(\angle)=\R(\angle^\perp)\cup \{ R_1(\angle)\} \cup \{ R_2(\angle)\}
$$
where $\R(\angle^\perp)$ contains all rays orthogonal to the plane of
the angle $\angle$, the rays $R_1(\angle )$ and $R_2(\angle )$
are orthogonal to the
first and second side of the oriented angle $\angle$ respectively.
Evidently, the set $\R(\angle)$ and the ordered pair
of rays $(R_1(\angle ), R_2(\angle ))$ define
the oriented angle $\angle$ uniquely.
We define the weight $\sigma (\angle )$ by the formula:
$$
\sigma (\angle)=
\cases
2/3, &\text{if $1\le \rho (R_1(\angle ), R_2(\angle ))\le d$,}\\
1/2, &\text{if $d+1\le \rho (R_1(\angle ), R_2(\angle ))\le 2d+1$,}\\
0, &\text{if $2d+2\le \rho (R_1(\angle ), R_2(\angle ))$.}
\endcases
$$
Here we take the distance in the graph $G(\R(\angle))$.
Let us prove conditions of the Lemma 1.4 with the constants
$C=(2/3)C_1+C_2/2$ and $D=0$.

The condition (1) follows from the condition (b) of the Theorem. We remark
that rays $R_1(\angle), R_2(\angle)$ do not belong to the set
$R(\gamma^\perp)$.

Let us prove the condition (2).

Let $\gamma _3$ be a 2-dimensional triangle face (triangle) of
$\gamma $. The set ${\Cal R }(\gamma _3)$ of all extremal rays orthogonal
to points of $\gamma_3$ is the union of the set
$\R(\gamma _3^\perp)$ of extemal rays, which
are orthogonal to the plane of the triangle
$\gamma _3$,
and rays $R_1, R_2, R_3$, which are orthogonal to sides of
the triangle $\gamma _3$. Union of the set
$\R(\gamma _3^\perp)$ with any two rays from
 $R_1, R_2, R_3$  is extremal, since it is orthogonal to
a vertex of $\gamma _3$.
On the other hand, the set $\R(\gamma _3)=
\R(\gamma _3^\perp)\cup \lbrace R_1, R_2, R_3 \rbrace$
is not extremal, since it is not orthogonal to a point of $\M(X)$.
Indeed, the set of all points of $\M(X)$, which are orthogonal
 to the set
$\R(\gamma _3^\perp)\cup \lbrace R_2, R_3 \rbrace$,
$\R(\gamma _3^\perp)\cup \lbrace R_1, R_3 \rbrace$,
or
$\R(\gamma _3^\perp)\cup \lbrace R_1, R_2\rbrace$ is the
vertex $A_1, A_2$ or $A_3$ respectively of the triangle
$\gamma _3$, and
the intersection of these sets of vertices is empty.
Thus, there exists an $E$-set $\La \subset \R(\gamma _3)$,
which contains the set of
rays $\{ R_1, R_2, R_3\} $. By the condition (a), the graph
$G(\La )$
contains a shortest oriented path $s$ of the
length $\le d$ which connects the rays
$R_1, R_3$. If this path does not contain the ray $R_2$, then the
 oriented angle of $\gamma _3$ defined by the set
$\R(\gamma _3^\perp )\cup \lbrace R_1, R_3 \rbrace$ and
the pair $(R_1, R_3)$ has the weight 2/3. If this path contains the
ray $R_2$, then the oriented angle of $\gamma _3$ defined by the
set
$\R(\gamma _3^\perp )\cup \lbrace R_1, R_2 \rbrace$ and
the pair $(R_1, R_2)$ has the weight 2/3. Thus, we proved that the side
$A_2A_3$ of the triangle $\gamma _3$ defines an oriented angle of
 the
triangle with the weight 2/3 and the first
side $A_2A_3$ of the oriented angle. The triangle has three sides.
It follows the condition (2) of the Lemma 1.4 for the triangle.

Let $\gamma _4$ be a 2-dimensional quadrangle face (quadrangle)
of $\gamma$. In this case,
$$
\R(\gamma _4)=
\R(\gamma _4^\perp)\cup \lbrace R_1, R_2, R_3, R_4
\rbrace
$$
where $\R(\gamma _4^\perp)$ is the set of all
extremal rays which are orthogonal to the plane of the quadrangle
and the rays $R_1, R_2, R_3, R_4$ are orthogonal to
the consecutive sides of the quadrangle. As above, one can see that the
sets
$\R(\gamma _4^\perp )\cup \lbrace R_1, R_3\rbrace$,
$\R(\gamma _4^\perp )\cup \lbrace R_2, R_4\rbrace$
are not extremal, but the sets
$\R(\gamma _4^\perp )\cup \lbrace R_1,R_2\rbrace$,
$\R(\gamma _4^\perp )\cup \lbrace R_2,R_3\rbrace$,
$\R(\gamma _4^\perp )\cup \lbrace R_3,R_4\rbrace$,
and
$\R(\gamma _4^\perp ) \cup \lbrace R_4,R_1\rbrace$
are extremal. It follows that there are $E$-sets $\La, {\Cal N}$
such that
$\lbrace R_1, R_3 \rbrace \subset \La \subset
\R(\gamma _4^\perp )\cup \lbrace R_1, R_3\rbrace$
and
$\lbrace R_2, R_4 \rbrace \subset {\Cal N}\subset
\R(\gamma _4^\perp)\cup \lbrace R_2, R_4\rbrace$.
By the Lemma 1.1, there
exist rays $R\in \La$ and $Q\in {\Cal N}$ such  that $RQ$ is
an arrow. By
the condition (a) of the Theorem, one of the rays $R_1, R_3$ is
joined by an oriented path $s_1$ of the length $\le d$ with the ray
$R$ and this path does not contain another ray from $R_1,R_3$
(here $R$ is the terminal of the path $s_1$). We can suppose that this
 ray is $R_1$ (otherwise, one
 should replace the ray $R_1$ by the ray $R_3$).  As above,
we can suppose that the ray $Q$ is
connected by the oriented path $s_2$ of the length $\le d$ with the
ray  $R_2$ and this path does not contain the ray $R_4$. The path
 $s_1RQs_2$ is an oriented path of the length $\le 2d+1$ in the
oriented graph
$G(\R(\gamma _4^\perp))\cup \lbrace R_1,R_2\rbrace)$.
It follows that the oriented angle of the
quadrangle $\gamma _4$, such that
consecutive sides of this angle are orthogonal to the rays
$R_1$ and $R_2$ respectively,
has the weight $\ge 1/2$. Thus, we proved that for a pair
of opposite sides of $\gamma_4$ there exists an oriented angle with
weight
$\ge 1/2$ such that the first side of this oriented angle is one
of this
opposite sides of the quadrangle. A quadrangle has two
pairs of opposite sides. It follows that the sum of weights of
oriented angles of $\gamma _4$ is $\ge 1$. It proves the condition (2)
of the Lemma 1.4 and the Theorem.
\enddemo

Below, we apply the Theorem 1.2 to 3-folds.

\newpage

\head
CHAPTER 2. THREEFOLDS
\endhead

\subhead
1. Contractible extremal rays
\endsubhead

We consider normal projective 3-folds $X$ with $\Bbb Q$-factorial
singularities.

Let $R$ be an extremal ray of Mori polyhedron $\overline {NE}(X)$
 of
$X$. A morphism $f:X\rightarrow Y$ on a normal projective
variety $Y$ is called the {\it contraction} of the ray $R$ if for
 an
irreducible curve $C$ of $X$ the
image $f(C)$ is a point iff $C\in R$. The contraction $f$
is defined by a linear system $H$ on $X$ ($H$ gives
 the nef element of $N^1(X)$, which we denote by $H$ also). It
follows that an irreducible curve $C$ is contracted iff
$C\cdot H=0$. We assume that the contraction $f$ has properties:
$f_\ast {\Cal O}_X={\Cal O}_Y$ and the sequence
$$
0\rightarrow {\Bbb R}R\rightarrow N_1(X)\rightarrow N_1(Y)
\rightarrow 0
\tag1--1
$$
is exact where the arrow $N_1(X)\rightarrow N_1(Y)$ is $f_\ast$.
An extremal ray $R$ is called {\it contractible} if there exists
its contraction $f$ with these properties.

The number $\kappa (R)=\dim Y$ is called {\it Kodaira dimension}
of the contractible extremal ray $R$.

A face $\gamma $ of $\overline {NE}(X)$ is called {\it contractible}
 if there exists a morphism $f:X\rightarrow Y$ on a normal projective
variety $Y$ such that $f_\ast \gamma =0$,
$f_\ast {\Cal O}_X={\Cal O}_Y$ and $f$  contracts  curves
lying in $\gamma $ only. The $\kappa (\gamma )=\dim Y$ is called
{\it Kodaira dimension of} $\gamma $.

 Let $H$ be a general nef element orthogonal to a face $\gamma $
of Mori polyhedron. {\it Numerical Kodaira dimension of} $\gamma$ is
defined by the formula
$$
\kappa_{num}(\gamma )=
\cases
3, &\text{if $H^3>0$;}\\
2, &\text{if $H^3=0$ and $H^2\not\equiv 0$;}\\
1, &\text{if $H^2\equiv 0$.}
\endcases
$$
It is obvious that for a contractible face $\gamma $ we have
$\kappa_{num}(\gamma )\ge \kappa (\gamma )$. In particular,
$\kappa_{num}(\gamma )=\kappa (\gamma )$ for a contractible face
$\gamma $ of Kodaira dimension $\kappa (\gamma )=3$.

\vskip 8pt

\subhead
2. Pairs of extremal rays of Kodaira dimension three
lying in contractible faces of $\overline {NE}(X)$
 of Kodaira dimension three
\endsubhead

Further $X$ is a projective normal threefold with $\Bbb Q$-factorial
singularities.

\proclaim{Lemma 2.2.1} Let $R$ be a contractible extremal ray of
Kodaira dimension 3 and $f:X\rightarrow Y$ its contraction.

Then there are three possibilities:

(I) All curves $C\in R$ fill an irreducible Weil divisor $D(R)$,
the contraction $f$ contracts $D(R)$ in a point and $R\cdot D(R)<0$.

(II) All curves $C\in R$ fill
an irreducible Weil divisor $D(R)$, the
contraction $f$ contracts $D(R)$ on an irreducible curve and
$R\cdot D(R)<0$.

(III) (small extremal ray)
All curves $C\in R$ give a finite set of irreducible curves
and the contraction $f$ contracts these curves in points.
\endproclaim

\demo{Proof} Assume that some curves of $R$ fill
an irreducible divisor $D$. Then $R\cdot D<0$
(this inequality follows from the Proposition 2.2.6 below). Suppose
that $C\in R$ and $D$ does not contain $C$. It follows that
$R\cdot D \ge 0$. We get a contradiction. It follows the Lemma.
\enddemo

According to the Lemma 2.2.1, we say that an extremal ray $R$ has the
{\it type (I), (II) or (III) (small)} if it is contractible of Kodaira
dimension 3 and the statements (I), (II) or (III) respectively hold.

\proclaim{Lemma 2.2.2}
Let $R_1$ and $R_2$ are two different extremal
rays of the type (I).
Then divisors $D(R_1)$ and $D(R_2)$ do not intersect one another.
 \endproclaim

\demo{Proof} Otherwise, $D(R_1)$ and $D(R_2)$ have a common curve and
the rays $R_1$ and $R_2$ are not different.
\enddemo

For an irreducible Weil divisor $D$ on $X$ let
$$
\nem (X,D)=(\hbox{image}\ \nem (D)) \subset \nem (X).
$$

\proclaim{Lemma 2.2.3}
Let $R$ be an extremal ray of the type (II),
and $f$ its contraction.

Then $\overline {NE}(X,D(R))=R + {\Bbb R}^+S$, where
${\Bbb R}^+ f_* S={\Bbb R}^+(f(D)).$
\endproclaim

\demo{Proof}
This follows at once from the exact sequence (1.1).
\enddemo

\proclaim{Lemma 2.2.4}
Let $R_1$ and $R_2$ are two different
 extremal rays of the type (II) such that the divisors $D(R_1)=D(R_2)$.

Then for
$D=D(R_1)=D(R_2)$ we have:
$$\overline {NE}(X, D)=R_1+R_2.$$
In particular, do not exist three different extremal rays of
the type (II) such that their divisors are coincided.
\endproclaim

\demo{Proof}
This follows from the Lemma 2.2.3.
\enddemo

\proclaim{Lemma 2.2.5}
Let $R$ be an extremal ray
of the type (II) and $f$ its contraction.

Then there does not exist more than one extremal ray $Q$ of
the type (I) such that $D(R)\cap D(Q)$ is not empty. If $Q$ is this
ray, then $D(R)\cap D(Q)$ is a curve and any irreducible component
of this curve is not contained in fibers of $f$.
\endproclaim

\demo{Proof} The last statement is obvious. Let us proof the first
one.
Suppose that  $Q_1$ and $Q_2$ are two different extremal
rays of the type (I) such that $D(Q_1)\cap D(R)$ and
$D(Q_2)\cap D(R)$ are not empty. Then the
plane angle $\overline {NE}(X,D(R))$ (see the Lemma 2.2.3)
contains three different extremal rays: $Q_1, Q_2$ and $R$. It is
impossible.
\enddemo

The following key Proposition is very important.

\proclaim{Proposition 2.2.6}
Let $X$ be a projective 3-fold with
$\Bbb Q$-factorial singularities, $D_1,...,D_m$ irreducible
divisors on $X$ and $f:X\rightarrow Y$ a surjective morphism such
 that
$\dim X=\dim Y$ and $\dim f(D_i)<\dim D_i$.
 Let $y\in f(D_1)\cap ...\cap f(D_m)$.

Then there are $a_1>0,\ ...,\ a_m>0$ and an open $U$, $y\in U\subset
f(D_1)\cup ...\cup f(D_m)$, such that
$$
C\cdot (a_1D_1+...+a_mD_m)<0 $$
if a curve $C\subset D_1\cup ...\cup D_m$ belongs to a non-trivial
algebraic family of curves on $D_1\cup ...\cup D_m$
and $f(C)=point\  \in U$.
\endproclaim

\demo{Proof} It is the same as for the well-known case
of surfaces (but, for surfaces, it is not necessary to suppose that
$C$ belongs to a nontrivial algebraic family). Let $H$ be an irreducible
ample divisor on $X$ and
$H'=f_\ast H$. Since $\dim f(D_i)<\dim D_i$, it follows that
 $f(D_1)\cup ...\cup f(D_m)\subset H'$. Let $\phi $ be a
non-zero rational
function on $Y$ which is regular in a neighbourhood $U$ of $y$ on
$Y$ and is equal to zero on the divisor $H'$. In the open set
$f^{-1}(U)$ the divisor
$$
(f^\ast \phi )=\sum_{i=1}^m a_iD_i+\sum_{j=1}^n b_jZ_j
$$
where all $a_i>0$ and all $b_j>0$. Here every divisor $Z_j$ is
different from any divisor $D_i$. We have
$$0=C\cdot \sum_{i=1}^m a_iD_i+C\cdot \sum_{j=1}^n b_jZ_j.$$
Here $C\cdot (\sum_{j=1}^n b_jZ_j)>0$ since C belongs to a nontrivial
algebraic family of curves on a surface
$D_1\cup ...\cup D_m$ and one
of the divisors $Z_j$ is the hyperplane section $H$.
\enddemo

\proclaim{Lemma 2.2.7}
Let $R_1, R_2$ are two extremal
rays of the type (II), divisors $D(R_1)$,
\linebreak
$D(R_2)$ are different
and $D(R_1)\cap D(R_2)\not= \emptyset$. Assume that $R_1,R_2$ belong
to a contractible face of $\overline {NE}(X)$ of Kodaira dimension 3.
 Let $0\not= F_1\in R_1$ and $0\not= F_2\in R_2$.

Then
$$
(F_1\cdot D(R_2))(F_2\cdot D(R_1))<(F_1\cdot D(R_1))(F_2\cdot
D(R_2)).
$$
\endproclaim

\demo{Proof}
Let $f$ be the contraction of a face of  Kodaira dimension 3,
which contains both rays $R_1,R_2$. By the proposition 2.2.6,
there are $a_1>0, a_2>0$ such that
$$
a_1(F_1\cdot D(R_1))+a_2(F_1\cdot D(R_2))<0\ \text{and}\
a_1(F_2\cdot D(R_1))+a_2(F_2\cdot D(R_2))<0.
$$
Or
$$
-a_1(F_1\cdot D(R_1))>a_2(F_1\cdot D(R_2))\ \text{and}\
-a_2(F_2\cdot D(R_2))>a_1(F_2\cdot D(R_1))
$$
where $F_1\cdot D(R_1)<0, F_2\cdot D(R_2)<0$ and
$F_1\cdot D(R_2)>0, F_2\cdot D(R_1)>0$.
Multiplying inequalities above, we obtain the Lemma.
\enddemo

\subhead
3. A classification of extremal sets of extremal rays
which contain extremal rays
of the type (I) and  simle extremal rays of the type (II)
\endsubhead

As above, we assume that $X$ is a projective normal
3-fold with $\Bbb Q$-factorial singularities.

\definition{Definition 2.3.1}
An extremal ray $R$ of the type (II) is
called {\it simple} if
$$
R\cdot (D(R)+D)\ge 0
$$
for any irreducible divisor $D$ such that $R\cdot D>0.$
\enddefinition

The following statement gives a simple sufficient condition for an
extremal ray to be simple.

\proclaim{Proposition 2.3.2}
Let $R$ be an extremal ray of
the type (II) and $f:X\rightarrow Y$ the contraction of $R$.
Suppose that the curve
$f(D(R))$ is not contained in the set of singularities of $Y$.

Then

(1) the ray $R$ is simple;

(2) if $X$ has only
isolated singularities, then a general element
$C$ of the ray $R$ (a general fiber of the morphism $f\mid D(R)$)
is isomorphic to $\Bbb P^1$ and the divisor $D(R)$ is non-singular
along $C$. If additionally $R\cdot K_X<0$, then
$C\cdot D(R)=C\cdot K_X=-1$.

(3) In particular, both statements (1) and (2) are true if $X$ has
terminal singularities and $R\cdot K_X<0$.
\endproclaim

\demo{Proof}
Let $D$ be an irreducible divisor  on $X$ such that
$R\cdot D>0$. Since $R\cdot D(R)<0$, the divisor $D$ is
different from $D(R)$ and the intersection $D\cap D(R)$ is
a curve which does not belong to $R$. Then $D'=f_\ast (D)$
is an irreducible divisor on $Y$ and $\Gamma =f(D(R))$ is
a curve on $D'$. Let $y\in \Gamma $ be a non-singular point of $Y$.
Then the divisor $D'$ is
defined by some local equation $\phi $ in a neighbourhood $U$ of $y$.
Evidently, in the open set $f^{-1}(U)$ the divisor
$$
(f^\ast \phi )=D+m(D(R))
$$
where the integer $m\ge 1$. Let a curve $C\in R$ and
$f(C)=y\in U\cap f(D(R))$. Then
$0=C\cdot (D+m(D(R)))=C\cdot (D+D(R))+C\cdot (m-1)(D(R))$. Since
$m\ge 1$ and $C\cdot D(R)<0$, it follows that $C\cdot (D+D(R))\ge
 0$.

Let us prove (2). Let us consider a linear system $\mid H\mid $ of hyperplane
sections on $Y$ and the corresponding linear systems
on the resolutions of singularities of $Y$ and $X$. Let us apply
Bertini's theorem (see, for example, \cite{H, ch. III, Corollary 10.9
and the Exercise 11.3})
to this linear systems. Singularities of $X$ and $Y$ are isolated.
Then by Bertini theorem,  for a general element $H$ of $\mid H\mid $
we obtain that
(a) $H$ and $f^{-1}(H)$ are irreducible and non-singular;
(b) $H$ intersects $\Gamma $ transversely in non-singular points
of $\Gamma $. Let us consider the corresponding birational
morphism $f'=f\mid H'\ :\ H'\rightarrow H$ of
the non-singular irreducible surfaces. It is a composition of
blowing ups in non-singular points. Thus, fibers of $f'$ over
$H\cap \Gamma $ are trees of non-singular rational curves.
The exceptional curve of the first of these blowing ups is
identified with the fiber of the projectivization of the normal
bundle
${\Bbb P}({\Cal N}_{\Gamma /Y})$. Thus, we obtain a rational map
over
the curve $\Gamma $
$$
\phi :{\Bbb P}({\Cal N}_{\Gamma /Y})\rightarrow D(R)
$$
of the irreducible surfaces. Evidently, it is the injection
in the general point of ${\Bbb P}({\Cal N}_{\Gamma /Y})$.
It follows that $\phi $ is a birational isomorphism of the surfaces.
Since $\phi $ is a birational map over the curve $\Gamma $,
it follows
that the general fibers of this maps are birationally isomorphic.
 It
follows that a general fiber of $f'$ is $C\simeq \Bbb P^1$. Since
$C$ is  non-singular and is an intersection of the non-singular
 surface $H'$ with the surface $D(R)$, and since
$X$ has only isolated singularities, it follows
that $D(R)$ is non-singular along the general curve $C$.

The $X$ and $D(R)$ are non-singular along $C\simeq {\Bbb P}^1$ and the
curve $C$ is non-singular.
Then the canonical class $K_C=(K_X+D(R))\mid C$ where
both divisors $K_X$ and $D(R)$ are Cartier
divisors on $X$ along $C$.
It follows that $-2=\text{deg}~K_C=K_X\cdot C+D(R)\cdot C$, where
 the
both numbers $K_X\cdot C$ and $D(R)\cdot C$ are negative integers.
Then $D(R)\cdot C=K_X\cdot C=-1.$

If $X$ has terminal sungularities and $R\cdot K_X <0$, then $Y$
has terminal singularities too (see, for example, \cite{Ka1}). Moreover,
3-dimensional terminal singularities are isolated.
 From (1), (2), the last statement of the Proposition follows.

In connection with Proposition 2.3.2, also see \cite{Mo2, 1.3 and
2.3.2} and \cite{I, Lemma 1}.
\enddemo

Let $R_1, R_2$ are two extremal rays of the type (I) or (II). They are
joined if $D(R_1)\cap D(R_2)\not= \emptyset$. It defines
{\it connected components} of a set of extremal rays of the type
(I) or (II).

We recall (see Chapter I) that a set $\E$ of extremal rays is called
{\it extremal}  if it is contained in a face of $\nem (X)$. We say that
$\E$ is {\it extremal of Kodaira dimension 3} if it is contained in
a face of numerical Kodaira dimension $3$ of $\nem(X)$.

We prove the following classification result.

\proclaim{Theorem 2.3.3} Let $\E=\{ R_1, R_2,..., R_n\} $
be an extremal set of extremal rays of the type (I) or (II). Suppose that
every extremal ray of $\E$ of the type (II) is simple. Assume that $\E$
is contained in a contractible face with Kodaira dimention $3$ of
$\nem (X)$. (Thus, $\E$ is extremal of Kodaira dimention $3$.)

Then every connected component of $\E$ has a type
$\Aa_1, \bb_2, \cc _m$ or $\dd_2$ below (see figure 1).

($\Aa_1$) One extremal ray of the type (I).

($\bb_2$) Two different extremal rays $S_1, S_2$ of the type (II) such that
their divisors $D(S_1)=D(S_2)$ are coincided.

($\cc_m$) $m\ge 1$ extremal rays $S_1, S_2,..., S_m$ of the type (II)
such that their divisors
$D(S_2), D(S_3),..., D(S_m)$ do not intersect one another, and
$S_1\cdot D(S_i)=0$ and $S_i\cdot D(S_1)>0$ for
$i=2,...,m$.

($\dd_2$) Two extremal rays $S_1, S_2$, where $S_1$ is of the type (II)
and $S_2$ of the type (I), $S_1\cdot D(S_2)>0$ and
$S_2\cdot D(S_1)>0$.
Either $S_1\cdot (b_1D(S_1)+b_2D(S_2))<0$ or
$S_2\cdot (b_1D(S_1)+b_2D(S_2))<0$ for any $b_1, b_2$ such that
$b_1\ge 0, b_2\ge 0$ and one of $b_1, b_2$ is
not zero.

The following inverse statement is true: If
$\E=\{ R_1, R_2,...,R_n\}$ is a connected set of extremal rays of the
type (I) or (II) and $\E$ has a type
$\Aa_1, \bb_2, \cc _m$ or $\dd_2$ above,
then $\E$ generates a simplicial face
$R_1+...+R_n$ of the dimension $n$ and numerical Kodaira dimension $3$ of
$\nem (X)$. In particular, extremal rays of the set $\E$ are linearly
independent.
\endproclaim

\demo{Proof} Let us prove the first statement. We can suppose that
$\E$ is connected. We have to prove that then $\E$ has the type
$\Aa_1, \bb_2, \cc _m$ or $\dd_2$. If $n=1$, this is obvious.

Let $n=2$. From Lemma 2.2.2, it follows that one of the rays $R_1, R_2$
has the type (II). Let $R_1$ has the type (II) and $R_2$ the type (I).
Since $D(R_1)\cap D(R_2)\not= \emptyset$, then evidently
$R_2\cdot D(R_1)>0$. If $R_1\cdot D(R_2)=0$, then the curve
$D(R_1)\cap D(R_2)$ belongs to the ray $R_1$. It follows that the rays
$R_1$ and $R_2$ contain the same curve. We get a contradiction. Thus,
$R_1\cdot D(R_2)>0$. The rays $R_1, R_2$ belong to a contractible face
of Kodaira dimension $3$ of Mori polyhedron. Let $f$ be a contraction of
this face. By the Lemma 2.2.3, $f$ contracts the divisors
$D(R_1), D(R_2)$ in a same point. By Proposition 2.2.6, there exist
positive $a_1, a_2$ such that $R_1\cdot (a_1D(R_1)+a_2D(R_2))<0$ and
$R_2\cdot (a_1D(R_1)+a_2D(R_2))<0$. Now suppose that for some $b_1>0$ and
$b_2>0$ the inequalities
$R_1\cdot (b_1D(R_1)+b_2D(R_2))\ge 0$ and
$R_2\cdot (b_1D(R_1)+b_2D(R_2))\ge 0$ hold.
There exists $\lambda >0$ such that $\lambda b_1\le a_1,
\lambda b_2 \le a_2$ and one of these inequalities is an equality. For
example, let $\lambda b_1=a_1$. Then
$$
R_1\cdot (a_1D(R_1)+a_2D(R_2))
=R_1\cdot \lambda (b_1D(R_1)+b_2D(R_2))+
R_1\cdot (a_2-\lambda b_2)D(R_2)\ge 0.
$$
We get a contradiction. It proves that in this case
$\E$ has the type $\dd_2$.

\midinsert
\vskip 2cm
\includegraphics{fanofig1.ps}
\vskip 3cm
\endinsert

 \centerline{Figure 1.}

\vskip5pt

Now assume that both rays $R_1, R_2$ have the type (II). Since the rays
$R_1, R_2$ are simple, from Lemma 2.2.7, it follows that either
$R_1\cdot D(R_2)=0$ or the $R_2\cdot D(R_1)=0$. If both these equalities
hold, the rays $R_1, R_2$ have a common curve. We get a contradiction.
Thus, in this case, $\E$ has the type $\cc_2$.

Let $n=3$. Every proper subset of $\E$ has connected components of types
$\Aa_1, \bb_2, \cc _m$ or $\dd_2$. Using Lemmas 2.2.2 ---2.2.5, one can see
very easily that either $\E$ has the type $\cc_3$ or we have the
following case:

The rays
$R_1, R_2, R_3$ have the type (II), every two element subset of $\E$ has the
type $\cc_2$ and we can find a numeration such that
$ R_1\cdot D(R_2)>0, R_2\cdot D(R_3)>0, R_3\cdot D(R_1)>0$. Let $f$ be
a contraction of the face $\gamma$. By Lemma 2.2.3, $f$ contracts the
divisors $D(R_1), D(R_2), D(R_3)$ in a one point. By Proposition 2.2.6,
there are positive $a_1, a_2, a_3$ such that
$$
R_i\cdot (a_1D(R_1)+a_2D(R_2)+a_3D(R_3))<0
$$
for $i=1,2,3$. On the other hand, from simplicity of the rays
$R_1, R_2, R_3$, it follows that
$$
R_i\cdot (D(R_1)+D(R_2)+D(R_3))\ge 0.
$$
Let $a_1=\min \{ a_1, a_2, a_3\}$. From the last inequality,
$$
R_1\cdot (a_1D(R_1)+a_2D(R_2)+a_3D(R_3))=
$$
$$
=R_1\cdot a_1(D(R_1)+D(R_2)+D(R_3))+R_1((a_2-a_1)D(R_2)+(a_3-a_1)D(R_3))\ge 0.
$$
We get a contradiction with the inequality above.

Let $n>3$. We have proven that every two or three element subset of $\E$ has
connected components of types
$\Aa_1, \bb_2, \cc _m$ or $\dd_2$. It follows very easily that then
$\E$ has the type $\cc_n$ (we suppose that $\E$ is connected).

Let us prove the inverse statement.

For the type $\Aa_1$ this is obvious.

Let $\E$ has the type $\bb_2$. Since the rays $S_1,S_2$ are extremal of
Kodaira dimension $3$, there are $nef$ elements $H_1, H_2$ such that
$H_1\cdot S_1=H_2\cdot S_2=0, {H_1}^3>0, {H_2}^3>0$. Let
$0\not=C_1\in S_1$ and $0\not=C_2\in S_2$. Let $D$ be a divisor of the rays
$S_1$ and $S_2$. Let us consider a map
$$
(H_1, H_2) \to H=
\tag3.1
$$
$$
=(-D\cdot C_2)(H_2\cdot C_1)H_1+(-D\cdot C_1)(H_1\cdot C_2)H_2+
(H_2\cdot C_1)(H_1\cdot C_2)D.
$$

For a fixed $H_1$, we get a linear map $H_2\to H$ of the set of $nef$
elements $H_2$ orthogonal to $S_2$ into the set of $nef$ elements $H$
orthogonal to $S_1$ and $S_2$. This map has a one dimensional kernel,
generated by
$(-D\cdot C_2)H_1+(H_1\cdot C_2)D$. It follows that
$S_1+S_2$ is a $2$-dimensional face of $\nem(X)$.

For a general $nef$ element $H=a_1H_1+a_2H_2+bD$ orthogonal to this face,
where $a_1, a_2, b>0$, we have
$H^3=(a_1H_1+a_2H_2+bD)^3\ge (a_1H_1+a_2H_2+bD)^2\cdot
(a_1H_1+a_2H_2)=(a_1H_1+a_2H_2+bD)\cdot (a_1H_1+a_2H_2+bD)\cdot
(a_1H_1+a_2H_2)\ge (a_1H_1+a_2H_2)^2\cdot (a_1H_1+a_2H_2+bD) \ge
(a_1H_1+a_2H_2)^3>0$, since
$a_1H_1+a_2H_2+bD$ and $a_1H_1+a_2H_2$ are $nef$. It follows that the face
$S_1+S_2$ has numerical Kodaira dimension $3$.

Let $\E$ has the type $\cc_m$. Let $H$ be a $nef$ element orthogonal to
the ray $S_1$. Let $0\not=C_i\in S_i$. Let us consider a map
$$
H \to H^\prime = H +
\sum_{i=2}^m {(-(H\cdot C_i)/(C_i\cdot D(S_i)))D(S_i)} .
\tag3.2
$$
It is a linear map of the set of $nef$ elements $H$ orthogonal to $S_1$
into the set of $nef$ elements $H^\prime$ orthogonal to the rays
$S_1, S_2,..., S_m$. The kernel of the map has the dimension $m-1$. It
follows that the rays $S_1,S_2,...,S_m$ belong to a face of
$\nem (X)$ of a dimension $\le m$. On the other hand, multiplying
rays $S_1,...,S_m$ on the divisors $D(S_1),...,D(S_m)$, one can see very
easily that the rays $S_1,...,S_m$ are linearly independent.
Thus, they generate a $m$-dimensional face of $\nem (X)$. Let us show that
this face is $S_1+S_2+...+S_m$.
To prove this, we show that every $m-1$ subset of $\E$ is contained in a face
of $\nem (X)$ of a dimension $\le m-1$.

If this subset contains the ray $S_1$, this subset has the type $\cc_{m-1}$.
By induction, we can suppose that
this subset belongs to a face of $\nem (X)$ of the
dimension $m-1$. Let us consider the subset $\{S_2, S_3,..., S_m\}$. Let
$H$ be an ample element of $X$. For the element $H$, the map (3.2) gives an
element $H^\prime$ which is orthogonal
to the rays $S_2,...,S_m$, but is not orthogonal
to the ray $S_1$. It follows that the
set $\{S_2,...,S_m\}$ belongs to a face of the Mori polyhedron of the
dimension $< m$. Like above, one can see that for a general $H$ orthogonal to
$S_1$ the element $H^\prime$ has $(H^\prime )^3 \ge H^3>0$.

Let $\E$ has the type $\dd_2$. Let $H$ be a $nef$ element orthogonal to the
ray $S_2$. Let $0\not= C_i\in S_i$. Let us consider a map
$$
H\to H^\prime = H+
\frac
{(H\cdot C_1)((-D(S_2)\cdot C_2)D(S_1)+(D(S_1)\cdot C_2)D(S_2))}
{(D(S_2)\cdot C_2)(D(S_1)\cdot C_1)-(D(S_1)\cdot C_2)(D(S_2)\cdot C_1)}.
\tag3.3
$$
Evidently,
$C_2\cdot ((-D(S_2)\cdot C_2)D(S_1)+(D(S_1)\cdot C_2)D(S_2))=0$.
 From this equality and the inequality from the definition of the system
$\dd_2$, it follows that $C_1\cdot ((-D(S_2)\cdot C_2)D(S_1)+
(D(S_1)\cdot C_2)D(S_2))<0$. Thus, the denominator from the formula (3.3)
is positive. Then (3.3) is a linear map of the set of
$nef$ elements $H$ orthogonal to the ray $S_2$ into the set of $nef$ elements
$H^\prime$ orthogonal to the rays $S_1, S_2$. Evidently, the map has a
one dimensional kernel. Thus, the rays $S_1$ and $S_2$ generate a two
dimensional face $S_1+S_2$ of Mori polyhedron. As above, for a general
element $H$ orthogonal to $S_2$ we have $(H^\prime)^3\ge H^3>0$.
\enddemo

\proclaim{Corollary 2.3.4}
Let ${\Cal E}=\{ R_1,R_2,...,R_n\}$
be an extremal set of extremal rays of the type (I) or (II) and every
extremal ray of ${\Cal E}$ of the type (II) is simple.
Assume that ${\Cal E}$ is contained in a contractible face
with Kodaira dimension 3 of the $\overline {NE}(X)$.
Let $m_1\ge 0, m_2\ge 0,...,m_n\ge 0$ and
at least one of $m_1,...,m_n$ is positive.

Then  there exists $i$, $1\le i\le n$, such that
$$
R_i\cdot (m_1D(R_1)+...+m_nD(R_n))<0.
$$
Thus, the condition (ii) from Chapter I is valid.
\endproclaim

\demo{Proof}
It is sufficient to prove this statement for the connected
$\Cal E$. For every type
${\goth A}_1$, ${\goth B}_2$, ${\goth C}_m$ and ${\goth D}_2$
of the Theorem 2.3.3,
one can prove it very easily.
\enddemo

Unfortunately, in general, the inverse statement of the
Theorem 2.3.3 holds only for connected extremal sets $\Cal E$.
We will give two cases when it is true for a non-connected $\Cal E$.

\definition{Definition 2.3.5}
A threefold $X$ is called
{\it strongly projective}
(respectively
{\it very strongly projective})
if the following statement holds: a set $\{ Q_1,...,Q_n\}$
of extremal rays of the type (II) is extremal of
Kodaira dimension 3
(respectively generates the simplicial face $Q_1+...+Q_n$ of
$\overline {NE}(X)$ of the dimension $n$ and Kodaira dimension 3)
if its divisors
$D(Q_1),...,D(Q_n)$ do not intersect one another.
\enddefinition

\proclaim{Theorem 2.3.6}
Let ${\Cal E}=\{ R_1,R_2,...,R_n\} $be a
 set of extremal rays of the type (I) or (II) such that every
connected component of $\Cal E$ has the type
${\goth A}_1, {\goth B}_2, {\goth C}_m$ or ${\goth D}_2$.

Then:

(1) $\Cal E$ is extremal of numerical Kodaira dimension 3
if and only if the same is true for
any subset of $\Cal E$ containing only extremal rays of the type (II)
whose divisors do not intersect one another. In particular,
it holds if $X$ is strongly projective.

(2) $\Cal E$ generates a simplicial face $R_1+...+R_n$ with
numerical Kodaira
dimension 3 of the Mori polyhedron if and only if the same is true for
any subset of $\Cal E$ containing only extremal rays of the type
(II) whose divisors do not intersect one another. In particular,
it is true if $X$ is very strongly projective.
\endproclaim

\demo{Proof} Let us prove (1). Only the inverse statement is non-trivial.
We prove it using an induction by $n$.
For $n=1$, the statement is obvious.

Assume that some connected component of $\Cal E$ has the type
${\goth A}_1$. Suppose that this component contains the ray $R_1$.
By induction, there exists a nef element $H$ such that $H^3>0$ and
$H\cdot R_i=0$  if $i>1$. Then there exists $k\ge 0$,
such that $H'=H+kD(R_1)$ is nef and $H'\cdot {\Cal E}=0$.
As above, one can prove that $(H')^3\ge H^3>0$.

Assume that some connected component of $\Cal E$ has the
type ${\goth B}_2$. Suppose that this component contains the rays
$R_1,R_2$ and $D(R_1)=D(R_2)=D$. Then, by induction, there are
nef elements $H_1$ and $H_2$ such that ${H_1}^3>0, {H_2}^3>0$ and
$H_1\cdot \{ R_1,R_3,...,R_n\} =0, H_2\cdot \{ R_2,R_3,...,R_n\}
=0$.
As for the proof of the inverse statement of the Theorem 2.3.3
in the case ${\goth B}_2$, there are $k_1\ge 0, k_2\ge 0, k_3\ge 0$ such
that the element $H=k_1H_1+k_2H_2+k_3D$ is nef, $H\cdot {\Cal E}=0$ and
$H^3>0$.

Assume that some connected component of $\Cal E$ has the type
${\goth C}_m, m>1$. We use notation of the Theorem 2.3.3
for this connected component. Let this is $\{ S_1,S_2,...,S_m\} $.
By the induction,
there exists a nef element $H$ such that $H$ is orthogonal to
${\Cal E}-\{ S_2,...,S_m\} $ and $H^3>0$.
As for the proof of the inverse
statement of the Theorem 2.3.3 in the case ${\goth C}_m$, there are
$k_2\ge 0,..., k_m\ge 0$ such that $H'=H+k_2D(S_2)+...+k_mD(S_m)$
 is nef,
$H'\cdot {\Cal E}=0$ and $(H')^3\ge H^3>0$.

Assume that some connected component of $\Cal E$ has the type
${\goth D}_2$. We use notation of the Theorem 2.3.3
for this connected component. Let this is $\{ S_1,S_2\}$.
By the induction, there exists nef element $H$ such
that $H^3>0$ and $H$ is orthogonal to ${\Cal E}-\{ S_1\} $.
As for the theorem 2.3.3, there are $k_1\ge 0, k_2\ge 0$ such that
$H'=H+k_1D(S_1)+k_2D(S_2)$ is nef, $H'\cdot {\Cal E}=0$ and
$(H')^3\ge H^3>0$.

If every connected component of $\Cal E$
has the type ${\goth C}_1$, then the statement holds by the
condition of the Theorem.

Let us prove (2). Only the inverse statement is non-trivial.
 We prove it
using an induction by $n$. For $n=1$ the statement is true. It is
sufficient to prove that $\Cal E$ is contained in a
face of a dimension $\le n$ of Mori polyhedron
because, by the induction, any its $n-1$ element subset generates
a simplicial face of the dimension $n-1$ of Mori polyhedron.

Assume that some connected component of $\Cal E$ has the type
${\goth A}_1$. Suppose that the ray $R_1$ belongs to this component
and $0\not= C_1\in R_1$. Let us consider the map
$$H\rightarrow H'=H+((H\cdot C_1)/(-D(R_1)\cdot C_1))D(R_1).$$
of the set of nef elements $H$ orthogonal to the set
$\{ R_2,...,R_n\} $ into the set of nef elements $H'$ orthogonal
to the
$\Cal E$. It is the linear map with one dimensional kernel.
Since, by the induction, the set $\{ R_2,...,R_n\} $ is contained
 in a
face of Mori polyhedron of the dimension $n-1$, it follows that $
\Cal E$
is contained in a face of the dimension $n$.

If $\Cal E$ has a connected component of the type
${\goth B}_2, {\goth C}_m, m>1,$ or ${\goth D}_2$,
a proof is the same if one uses the maps (3.1), (3.2) and (3.3) above.

If all connected components of $\Cal E$ have the type ${\goth C}_
1$,
the statement holds by the condition.
\enddemo

\remark{Remark 2.3.7}
Like the statement (1) of the Theorem 2.3.6, one can
prove that a set $\Cal E$ of extremal rays with connected
components of the type ${\goth A}_1,{\goth B}_2, {\goth C}_m$ or
${\goth D}_2$ is extremal if and only
if the same is true for any subset of $\Cal E$
containing only extremal rays of the type (II) whose divisors
do not intersect one another.
\endremark

The next statement is simple but important. To simplify notation, we
say that for a fixed $a_1,...,a_n$ we have a
{\it linear dipendence condition}
$$
a_1R_1+...+a_nR_n=0
$$
between extremal rays $R_1,...,R_n$ if there exist non-zero $C_i\in R_i$ such
that
$$
a_1C_1+...+a_nC_n=0.
$$

\proclaim{Proposition 2.3.8}
Assume that a set $\Cal E=
\{ R_1,R_2,...,R_m \} $ of extremal rays
has connected components of the type
${\goth A}_1,{\goth B}_2, {\goth C}_m$ or ${\goth D}_2$ and
there exists a linear dependence condition
$a_1R_1+a_2R_2+...+a_mR_m=0$ with
all $a_i\not=0$.

Then all connected components of $\E$ have the type $\bb _2$.
Let these components are ${\goth B}^1,...,{\goth B}^t$. Then
$t\ge 2$, and we can choose a
numeration such that ${\goth B}^i=\{ R_{i1},R_{i2}\} $ and the linear
dependence has a form
$$
a_{11}R_{11}+a_{21}R_{21}+...+a_{t1}R_{t1}=
a_{12}R_{12}+a_{22}R_{22}+...+a_{t2}R_{t2}.
$$
where all $a_{ij}>0$.
\endproclaim

\demo{Proof}
Let us multiply the
equality
$a_1R_1+a_2R_2+...+a_mR_m=0$
on divisors \linebreak $D(R_1),...,D(R_m)$. Then we get that
$a_{k}=0$ if the ray $R_k$ belongs to a connected component
of the type
${\goth A}_1, {\goth C}_m$ or ${\goth D}_2$. Thus, all connected
components of $\Cal E$ have the type $\bb_2$. Let these components are
$$
{\goth B}^1=\{ R_{11},R_{12}\},{\goth B}^2=\{ R_{21},R_{22}\}, \hdots ,
{\goth B}^t=\{ R_{t1},R_{t2}\}.
$$
Obviously, $t\ge 2$, and we can rewrite the linear dependence as
$$
a_{11}R_{11}+a_{12}R_{12}+a_{21}R_{21}+
a_{22}R_{22}+...+a_{t1}R_{t1}+a_{t2}R_{t2}=0,
$$
where all $a_{ij}\not= 0$. Multiplying this equation on all divisors
$D(R_{ij})$ and using inequalities
$R_{ij}\cdot D(R_{ij})<0$, we get the last statement of Proposition.
\enddemo

\subhead
4. A classification of $E$-sets of
extremal rays of the type (I) or (II)
\endsubhead

As above, we suppose that $X$ is a projective normal 3-fold with
$\Bbb Q$-factorial singularities.

We recall that a set $\Cal L$ of extremal rays is called
{\it $E$-set} if it is not extremal but any proper subset
 of $\Cal L$ is extremal (it is
contained in a face of $\overline {NE}(X)$).
Thus, an $E$-set is a minimal non-extremal set of extremal rays.

\proclaim{Theorem 2.4.1}
Let $\Cal L$ be a $E$-set of  extremal rays
of the type (I) or (II).
Suppose that every ray of the type (II)
of $\Cal L$ is simple and every proper subset of $\Cal L$ is
contained
in a contractible face of Kodaira dimension 3 of Mori polyhedron.

Then we have one of the following cases:

(a) $\Cal L$ is connected and ${\Cal L}=\{ R_1,R_2,R_3\} $, where
any $R_i$ has the type (II) and
every of 2-element subsets $\{ R_1,R_2\} ,
 \{ R_2,R_3\} ,\{ R_3,R_1\} $  of $\Cal L$ has the type
${\goth C}_2$. Here $R_1\cdot D(R_2)>0, R_2\cdot D(R_3)>0,
R_3\cdot D(R_1)>0$ but $R_2\cdot D(R_1)=R_3\cdot D(R_2)=
R_1\cdot D(R_3)=0$. The divisor
$D({\Cal L})=D(R_1)+D(R_2)+D(R_3)$ is nef.

(b) $\Cal L$ is connected and ${\Cal L}=\{ R_1,R_2\}$, where at least
one of the rays $R_1, R_2$ has the type (II). There are
positive $m_1,m_2$ such that $R\cdot  (m_1D(R_1)+m_2D(R_2))\ge 0$
for any extremal ray $R$ of the type (I) or simple
extremal ray of the type (II) on $X$. If the divisor $m_1D(R_1)+
m_2D(R_2)$ is not nef, both the extremal rays $R_1, R_2$ have the
type (II).

(c) $\Cal L$ is connected and ${\Cal L}=\{ R_1,R_2\} $ where
both $R_1$ and $R_2$
have the type (II) and  there exists the simple extremal
ray $S_1$ of the type (II) such that the rays $R_1, S_1$ define the
extremal set of the type ${\goth B}_2$ (it means that
$S_1\not= R_1$ but the divisors
$D(S_1)=D(R_1)$) and the rays $S_1, R_2$ define
the extremal set of the type ${\goth C}_2$, where
$S_1\cdot D(R_2)=0$ but $R_2\cdot D(S_1)>0$.  Here there do not exist
positive $m_1, m_2$ such that the divisor $m_1D(R_1)+m_2D(R_2)$ is nef,
since evidently $S_1\cdot (m_1D(R_1)+m_2D(R_2))<0$.
See figure 2 below.

(d) ${\Cal L}=\{ R_1,...,R_k\} $  where $k\ge 2$, all rays
$R_1, \hdots ,R_k$ have the type (II) and the divisors $D(R_1), \hdots ,
D(R_k)$
do not intersect one another. Any proper subset of ${\Cal L}$ is
contained in a contractible face of Kodaira dimension 3 of Mori
polyhedron but $\Cal L$ is not contained in a face of Mori polyhedron.
\endproclaim

\comment
\midspace{4cm} \caption{Figure 2}
\endcomment

\midinsert
\vskip 2cm
\includegraphics{fanofig2.ps}
\vskip 3cm
\endinsert

\centerline{Figure 2.}

\vskip5pt

\demo{Proof}
Let ${\Cal L}=\{ R_1,...,R_n\} $ be a $E$-set of extremal rays
satisfying the conditions of the Theorem. Let us consider two
cases.

The case 1. Let $\Cal L$ is not connected. Then  every connected
component of $\Cal L$ is extremal and, by the theorem 2.3.3,
it has the type ${\goth A}_1,{\goth B}_2, {\goth C}_m$ or ${\goth D}_
2$. If some of
these components does not have the type
${\goth C}_1$, then, by the statement
(1) of the Theorem 2.3.6, $\Cal L$ is extremal and we get the
contradiction. Thus, we get the case (d) of the Theorem.

The case 2. Let ${\Cal L}=\{ R_1,...,R_n\} $ is connected.

Let $n\ge 4$. By the Theorem 2.3.3, any proper subset of $\Cal L$
 has connected components of the type
${\goth A}_1,{\goth B}_2, {\goth C}_m$ or ${\goth D}_2$.
Like for the proof of the Theorem 2.3.3, it
follows that $\Cal L$ has the type
${\goth C}_n$.
By the Theorem 2.3.3, then $\Cal L$ is extremal. We get
the contradiction.

Let $n=3$. Then, like for the proof of the Theorem 2.3.3,
we get that $\Cal L$ has the type (a).

Let $n=2$ and ${\Cal L}=\{ R_1,R_2\} $. If both rays $R_1,R_2$
have the type (I), then, by the Lemma 2.2.2,
$\Cal L$ is not connected and we get the contradiction.

Let $R_1$ has the type (I) and $R_2$ has the type (II). Since the
set $\Cal L$ is not extremal, by the Theorem 2.3.3, there are
positive $m_1,m_2$ such that $R_1\cdot (m_1D(R_1)+m_2D(R_2))\ge 0
$ and
$R_2\cdot (m_1D(R_1)+m_2D(R_2))\ge 0$. By the Lemma 2.2.3, it follows
that $C\cdot (m_1D(R_1)+m_2D(R_2))\ge 0$ if the curve $C$ is
contained in the $D(R_1)\cup D(R_2)$. If $C$ is not contained in
$D(R_1)\cup D(R_2)$, then obviously $C\cdot (m_1D(R_1)+m_2D(R_2))
\ge 0$.
It follows, that the divisor $m_1D(R_1)+m_2D(R_2)$ is nef. Thus,
we get the case (b).

Let both rays $R_1, R_2$ have the type (II). If $D(R_1)=D(R_2)$, then
we get an extremal set $\{ R_1, R_2\}$ by the Theorem 2.3.3.
Thus, the divisors
$D(R_1)$ and $D(R_2)$ are different. By the Lemma 2.2.1, the curve
$D(R_1)\cap D(R_2)$ does not have
an irreducible component which belongs to both rays  $R_1$ and
$R_2$. Since rays $R_1, R_2$ are simple, it follows that
$R_1\cdot (D(R_1)+D(R_2))\ge 0$ and $R_2\cdot (D(R_1)+D(R_2))\ge 0$.
Let $R$ be an extremal ray of the type (I) or
simple extremal ray of the type (II). If the divisor
$D(R)$ does not coincide with the divisor $D(R_1)$ or $D(R_2)$,
then obviously $R\cdot (D(R_1)+D(R_2))\ge 0$. Thus, if there does not
exist
an extremal ray $R$ which has the same divisor
as the ray $R_1$ or $R_2$, we get the case (b).

Assume that $D(R)=D(R_1)$. Then, by the Lemma 2.2.5, the ray $R$
has
the type (II) too. If $R\cdot D(R_2)=0$, we get the case (c) of
the Theorem where $S_1=R$. If $R\cdot D(R_2) >0$, then
$R\cdot (D(R_1)+D(R_2))\ge 0$ since the ray $R$ is simple. Then
we get the case (b) of the Theorem.
\enddemo

\subhead
5. An application of the diagram method to Fano 3-folds with
terminal singularities
\endsubhead

We restrict ourselves by considering Fano 3-folds with $\Bbb Q$-factorial
terminal singularities, but it is possible
to formulate and prove corresponding results
for a negative part of Mori cone of
3-dimensional variety with $\Bbb Q$-factorial
terminal singularities like in \cite{N7} .

We recall that an algebraic 3-fold $X$ over $\Bbb C$ with $\Bbb Q$-factorial
singularities is called Fano if the anticanonical class $-K_X$ is ample.
By results of Kawamata \cite{Ka1} and Shokurov \cite{Sh}, any face of
$\nem (X)$ is contractible and $\nem (X)$ is generated by a finite set of
extremal rays if $X$ is a Fano 3-fold with terminal $\Bbb Q$-factorial
singularities.

\subsubhead
5.1. Preliminary results
\endsubsubhead

We need the following

\proclaim{Lemma 2.5.1}
Let $X$ be a Fano 3-fold with $\Bbb Q$-factorial terminal singularities.
Let
$\E=\{ R_1,...,R_n\}$
be a set of $n$ extremal rays on of the
type (II) and with disjoint divisors $D(R_1),...,D(R_n)$ on $X$.
(Thus, $\E$ has the type $n\cc_1$).

 If we suppose that the set $\E$ is not extremal,
then there exists a small extremal ray
$S$ and $i$, $1\le i \le n$, such that $S\cdot (-K_X+D(R_i))<0$ and
$S\cdot D(R_j)=0$
if $j\not=i$.

It follows that any curve of the ray $S$ belongs to the divisor $D(R_i)$.
\endproclaim

\demo{Proof} By Proposition 2.3.2, the divisor
$H=-K_X+D(R_1)+...D(R_n)$ is orthogonal to $\E$. Besides,
$H$ is $nef$ and $H^3>0$
if there does not exist a small
extremal ray $S$ with the property above. Then, $\E$ is extremal of
Kodaira dimension 3.
\enddemo

\definition{Definition 2.5.2} A set $\{ R, S \}$ of  extremal rays
has the type $\ee_2$ if the ray $R$ has type (II), the
extremal ray $S$
is small and $S\cdot D(R)<0$. (See Figure 3.)
\enddefinition

Thus, by Lemma 2.5.1, the set $R_1,...,R_n, S$ of extremal rays
contains a subset of the type $\ee_2$.

 From Proposition 2.3.2, any extremal ray of $X$ of the type (II) is
simple, and from Sections 3 and 4 we get a classification
of extremal sets and
$E$-sets of extremal rays of the type (I) and (II) on $X$.

\midinsert
\vskip 2cm
\includegraphics{fanofig3.ps}
\vskip 3cm
\endinsert

\centerline{Figure 3.}

\vskip5pt

We have the following general Theorem.

\proclaim{Theorem 2.5.3} Let $X$ be a Fano 3-fold with
$\Bbb Q$-factorial terminal singularities.
Let $\alpha $ be a face of $\nem (X)$.
Then we have the following possibilities:

(1) There exists a small extremal ray $S$ such that $\alpha +S$ is
contained in a face of $\nem (X)$ of Kodaira dimension 3.

(2) There are extremal rays $R_1, R_2$ of the type (II) and
a small extremal ray $S$ such that
$\alpha +R_1$ and $\alpha+R_2$ are contained in faces of $\nem (X)$ of
Kodaira dimension 3, the ray $R_2$ does not belong to $\alpha$, and one
of the sets $\{ R_1, S\}$ or $\{ R_2, S\}$ has the type $\ee_2$.

(3) The face $\alpha $ is contained in a face of $\nem (X)$ of
Kodaira dimension 1 or 2.

(4) There exists an E-set $\La =\{ R_1, R_2 \}$ such that
$R_1\not\subset \alpha$, $R_2\not\subset \alpha$, but
$\alpha +R_1$ and $\alpha +R_2$ are
 contained in faces of $\nem(X)$ of Kodaira dimension $3$.
The $\La $ satisfies the condition (c) of Theorem 2.4.1: Thus, both
extremal rays $R_1, R_2$ have the type (II) $R_1\cdot D(R_2)>0$ and
$R_2\cdot D(R_1)>0$ and there exists an
extremal ray $R_1^\prime$ of the type (II) such that
$D(R_1)=D(R^\prime _1)$ and $R^\prime_1 \cdot D(R_2)=0$.

(5) There are extremal rays $R_1,...,R_n$ of the type (II) such that
any of them
does not belong to $\alpha$,
$\alpha + R_1+...+R_n$ is contained in a face of $\nem (X)$ of
Kodaira dimension $3$ and
$$
\dim \alpha +R_1 +...+ R_n < \dim \alpha +n.
$$

(6) $\dim N_1(X)-\dim \alpha \le 12$.
\endproclaim

\demo{Proof} Let us consider the face $\gamma =\alpha ^\perp$ of
$\M (X)$ and apply Theorem 1.2 to this face $\gamma $. We have
$\dim \gamma =\dim N_1(X)-1-\dim \alpha$.

Assume that $\alpha $ does not satisfy the conditions
(1), (3) and (5). Then $\R(\gamma)$ contains extremal rays of the type
(I) or (II) only and $\M(X)$ is closed and simple in the face
$\gamma$. By  Proposition 2.3.2 and Theorem 2.3.3, any extremal subset
$\E$ of $\R(\gamma )$ has connected components of the types
$\Aa_1, \bb_2, \cc_n$ or $\dd_2$. By Corollary 2.3.4,
the condition (ii) is valid for extremal subsets of
$R(\gamma )$. Let $\La \subset \R(\gamma )$ be a $E$-set.
Assume that at least two elements $R_1, R_2\in \La$
don't belong to $\R (\gamma^\perp)$ and
for any proper subset $\La^\prime \subset \La$ we have that
$\La^\prime \cup \R (\gamma ^\perp)$ is extremal. Let us apply
Theorem 2.4.1 to $\La$.

Assume that $\La$ has the type (d). By  Lemma 2.5.1, one of extremal
rays $R_1$ of $\La$ together with some small extremal ray $S$ define
a set of the type $\ee_2$. Since $\{R_1\}\subset \La $ is a proper subset
of $\La$, the $\R(\gamma ^\perp )\cup \{R_1\}$ is extremal. Or
$\alpha +R_1$ is contained in a face of $\nem (X)$.
Since $\La$ has at least 2 elements which do not belong to
$\R(\gamma ^\perp)$, there exists another extremal ray  $R_2$ of $\La$ which
does not belong to $\R(\gamma ^\perp)$. Like above,
$\alpha + R_2$ is contained in a face of $\nem (X)$ of Kodaira dimension 3.
By definition of the case (d), both extremal rays $R_1, R_2$ have the
type (II).
Thus, we get the case (2) of Theorem.

Assume that $\La$ has the type (c). Then we get the case (4) of Theorem.

Assume that $\La =\{R_1,R_2\}$ has the type (b).
Suppose that the divisor $m_1D(R_1)+m_2(D(R_2)$
is not $nef$ (see the case (b) of Theorem 2.4.1). Then
there exists a small extremal ray $S$ such that
$S\cdot (m_1D(R_1)+m_2D(R_2))<0$.
It follows that one of the sets $\{ R_1,S\}$ or $\{R_2,S\}$ has the type
$\ee_2$. Thus, we get the case (2).

Assume that $\La=\{R_1, R_2, R_3\}$ has the type (a).
Then the divisor $D(R_1)+D(R_2)+D(R_3)$ is $nef$.

Thus, if we additionally exclude the cases (2) and (4), then all conditions of
the Theorem 1.2 are satisfied. By Theorems 2.4.1 and 2.3.3,
we can take $d=2$, $C_1=1$ and $C_2=0$. (See Figure 4 for graphs
$G(\E)$ corresponding to extremal sets $\E$ of the types
$\Aa_1, \bb_2, \cc_m$ and $\dd_2$.)
Thus,  by Theorem 1.2, $\dim \gamma <34/3$. It follows that
$\dim N_1(X)-\dim \alpha \le 12.$
\enddemo

\midinsert
\vskip 2cm
\includegraphics{fanofig4.ps}
\vskip 3cm
\endinsert

\centerline{Figure 4}

\vskip5pt

\subsubhead
5.2. General properties of
configurations of extremal rays of the type $\bb_2$
\endsubsubhead

Let $\{ R_{11}, R_{12}\}$ be a set of extremal rays of the type $\bb_2$.
By Theorem 2.3.3, they define a 2-dimensional face
$R_{11}+R_{12}$ of $\nem (X)$.
Let
$\{ R_{21}, R_{22}\}$ be another set of
extremal rays of the type $\bb_2$.  Since two different
2-dimensional faces of
$\nem (X)$ may have only a common extremal ray,
the divisors $D(R_{11})=D(R_{12})$ and $D(R_{21})=D(R_{22})$  don't
have a common point. There exists the maximal set
$\{ R_{11}, R_{12}\} ,\{ R_{21}, R_{22}\},..,\{ R_{n1}, R_{n2}\}$ of
pairs of extremal rays of the type $\bb_2$.

\proclaim{Lemma 2.5.4} Any $t$ pairs
$\{ R_{11}, R_{12}\} ,\{ R_{21}, R_{22}\},..,\{ R_{t1}, R_{t2}\}$
of extremal rays of the type $\bb_2$ generate a face
$$
\sum_{i=1}^{t} \sum_{j=1}^{2} {R_{ij}}
\subset \nem (X)\subset N_1(X)
$$
of the Kodaira dimension $3$ of $\nem(X)$.
\endproclaim

\demo{Proof} This face is orthogonal to the nef divisor
$H = -K_X + D(R_{11}) +... + D(R_{t1})$ with $H^3\ge (-K_X)^3 > 0$.
\enddemo

\proclaim{Lemma 2.5.5} In notation above, there exists a changing of
 numeration of pairs of extremal rays $R_{i1}, R_{i2}$ such that
$R_{11}+\cdots R_{t1}$ is a simplicial face of $\nem (X)$.
\endproclaim

\demo{Proof} For $t=1$, it is obvious. Let us suppose that
$\theta=R_{11}+\cdots R_{(t-1)1}$ is a simplicial face of the face
$$
\alpha_{t-1}=\sum_{i=1}^{t-1} \sum_{j=1}^{2} {R_{ij}}.
$$
The face
$$
\alpha_t=\sum_{i=1}^{t} \sum_{j=1}^{2} {R_{ij}}
$$
has $\alpha_{t-1}$ as its face and does not coincide with the
face $\alpha_{t-1}$. It follows that there exists a face $\beta$
of $\alpha_t$ of the dimension $t$ such that
$\beta\not\subset \alpha_{t-1}$ but $\theta \subset \beta$ is a face
of $\beta$. It follows that all extremal rays of $\beta$ are the
extremal rays $R_{11},...,R_{(t-1)1}$ and some of extremal rays
$R_{t1},R_{t2}$. Assume that both extremal rays
$R_{t1}, R_{t2}$ belong to $\beta$. Then the extremal rays
$R_{11},...R_{(t-1)1}, R_{t1}, R_{t2}$ are linearly dependent, since
$\dim \beta =t$. By Proposition 2.3.8, it is impossible. Thus, only
one of extremal rays $R_{t1}, R_{t2}$ belongs to the face $\beta$.
Suppose that this is $R_{t1}$. Then
$\beta=R_{11}+\cdots R_{(t-1)1}+R_{t1}$ will be the face we were looking for.
\enddemo

We divide the maximal set
$\{ R_{11}, R_{12}\} ,\{ R_{21}, R_{22}\},..,\{ R_{n1}, R_{n2}\}$ of
pairs of extremal rays of the type $\bb_2$ on two parts:
$$
\{ R_{11}, R_{12}\} ,\{ R_{21}, R_{22}\},..,\{ R_{m1}, R_{m2}\}
$$
and
$$
\{ R_{(m+1)1}, R_{(m+1)2}\} ,\{ R_{(m+2)1}, R_{(m+2)2}\},..,
\{ R_{(m+k)1}, R_{(m+k)2}\}
$$
where $n=m+k$.
By definition, here the extremal rays $R_{i1}, R_{i2}$ belong to the
first part if and only if they are linearly independent from other extremal
rays from the set
$\{ R_{11}, R_{12}\} ,\{ R_{21}, R_{22}\},..,\{ R_{n1}, R_{n2}\}$.
Thus, extremal rays $R_{j1}, R_{j2}$
belong to the second part if they are linearly dependent from other extremal
rays from the set \linebreak
$\{ R_{11}, R_{12}\} ,\{ R_{21}, R_{22}\},..,\{ R_{n1}, R_{n2}\}$.

\proclaim{Lemma 2.5.6} Let $S$ be an extremal ray of the type (II)
such that $\{ R_{i1}, S\} $ define a configuration (c) of the
Theorem 2.4.1. Thus: $R_{i1}\cdot D(S)>0$, $S\cdot D(R_{i1})>0$ and
$R_{i2}\cdot D(S)=0$. Then
the extremal rays $R_{i1}, R_{i2}$ are linearly
independent from all other extremal
rays in
$\{ R_{11}, R_{12}\} ,\{ R_{21}, R_{22}\},..,\{ R_{n1}, R_{n2}\}$.
Thus, $1\le i \le m$. There does not exist a configuration of this type
with the ray $R_{i2}$. Thus, there does not exist an extremal ray
$S^\prime$ of the type (II) such that
$R_{i2}\cdot D(S^\prime)>0$, $S^\prime \cdot D(R_{i2})>0$ and
$R_{i1}\cdot D(S^\prime )=0$.
\endproclaim

\demo{Proof} The $R_{i1}+R_{i2}$ and $R_{i2}+S$ are 2-dimensional faces
of $\nem (X)$ with intersection by the extremal
ray $R_{i2}$. It follows that any curve of $D(S)$ belongs to the face
$R_{i2}+S$ (by Lemma 2.2.3). It follows that the divisor
$D(S)$ has no common point with the divisor $D(R_{j1})$ for any
other pair $R_{j1},R_{j2}$ for $j\not=i$.
Multiplying on $D(S)$ a linear relation of extremal rays
$R_{i1},R_{i2}$ with other extremal rays
$\{ R_{11}, R_{12}\} ,\{ R_{21}, R_{22}\},..,\{ R_{n1}, R_{n2}\}$
and using Proposition 2.3.8, we get that this linear relation does not
exist.

Let us suppose that there exists an extremal ray $S^\prime $ (see formulation
of the Lemma). Then $R_{i1}+S^\prime$ is another 2-dimensional face
of $\nem(X)$. Evidently, divisors $D(S)$ and $D(S^\prime)$ have a
non-empty intersection. Thus, faces $R_{i2}+S$ and $R_{i1}+S^\prime$
have a common ray. But it is possible only if $S=S^\prime$. Thus, we
get a contradiction, because $R_{i1}\cdot D(S)>0$ but
$R_{i1}\cdot D(S^\prime)=0$.
\enddemo

Using this Lemma 2.5.6, we can subdivide the first set
$$
\{ R_{11}, R_{12}\} ,\{ R_{21}, R_{22}\},..,\{ R_{m1}, R_{m2}\}.
$$
on sets
$$
\{ R_{11}, R_{12}\} ,\{ R_{21}, R_{22}\},..,\{ R_{m_11}, R_{m_12}\}
$$
and
$$
\{ R_{(m_1+1)1}, R_{(m_1+1)2}\},..,\{ R_{(m_1+m_2)1}, R_{(m_1+m_2)2}\}
$$
where $m_1+m_2=m$. Here $R_{i1},R_{i2}$ belong to the first part if  and
only if
there exists an extremal ray $S$ such that
$R_{i1}, S$ satisfies the condition of Lemma 2.5.6. By Lemma 2.5.6,
the numeration between extremal rays $R_{i1}$ and $R_{i2}$ is
then canonical.

Let us consider the second set
$$
\{ R_{(m+1)1}, R_{(m+1)2}\} ,\{ R_{(m+2)1}, R_{(m+2)2}\},..,
\{ R_{(m+k)1}, R_{(m+k)2}\}.
$$
We introduce an invariant
$$
\delta = \dim \sum_{i=m+1}^{m+k} \sum_{j=1}^{2} {R_{ij}} - k
$$
of $X$. Evidently,
$$
\text{either\ } k=\delta =0\ \text{or\ } k\ge 2\  \text{and\ } 1\le \delta<k.
$$
Thus,
$$
\dim \sum_{i=m+1}^{m+k} \sum_{j=1}^{2} {R_{ij}}=k+\delta .
$$
Let $$
\rho_0(X)=\dim N_1(X)-
\dim \sum_{i=1}^{n=m+k} \sum_{j=1}^{2} {R_{ij}}.
$$
Then
$$
\rho (X)=\dim N_1(X)=\rho_0(X)+2m+k+\delta.
$$
The invariants: $\rho_0(X), n, m, k, \delta, m_1, m_2$ are important
invariants of a Fano 3-fold $X$.

The following Lemma will be very useful:

\proclaim{Lemma  2.5.7}
Let $\E$ be an extremal set of extremal rays. Let
$$
\{ R_{11}, R_{12}\} \cup ...\cup \{ R_{t1}, R_{t2}\}
$$
be a set of different pairs of extremal rays of the type $\bb_2$.
Assume that $R\cdot D(R_{i1})=0$ for any $R\in \E$ and any
$i$, $1\le i \le t$.

Then there are extremal rays
$Q_1,...,Q_r$ such that the following conditions are valid:

(a) $r \le t$;

(b) For any $i$, $1\le i \le r$, there exists $j$, $1\le j \le t$,
such that $Q_i\cdot D(R_{j1})>0$ (in particular, $Q_i$ is different
from extremal rays of pairs of extremal rays
$\{R_{u1},R_{u2}\}$ of the type $\bb_2$);

(c) For any $j$, $1\le j\le t$, there exists an extremal
ray $Q_i$, $1\le i \le r$, such that
$$
Q_i \cdot D(R_{j1})>0;
$$

(d) The set
$\E\cup \{ Q_1,...,Q_r\}$
of extremal rays
is extremal, and extremal rays
$\{ Q_1,...,Q_r\}$ are linearly independent.
\endproclaim

\demo{Proof} Let us consider a linear subspace
$V\subset N_1(X)$ generated by all extremal rays
$\E$. Let us consider a factorisation map
$\pi : N_1(X)\to N_1(X)/V$. The cone $\pi (\nem (X))$ is generated
by images of extremal rays $T$ which together with all extremal rays
from $\E$ belong to faces of $\nem (X)$. There exists a curve $C$ on
$X$ such that $C\cdot D(R_{11})>0$. This curve $C$ (as
any element $x\in \nem (X)$) is a linear combination of extremal rays $T$
with non-negative coefficients and extremal rays from $\E$ with
real coefficients. We have $R\cdot D(R_{11})=0$ for any
extremal ray $R\in \E$. Thus,
there exists an extremal ray $T$ above
such that $T\cdot D(R_{11})>0$.
It follows that $T$ is different from extremal rays of pairs
of the type $\bb_2$. We take $Q_1=T$. By our construction, the set
$\E \cup \{ Q_1\}$ is extremal since it is contained in a face of
$\nem (X)$.
If $Q_1 \cdot D(R_{j1})>0$ for any $j$ such that
$1 \le j \le t$, then $r=1$, and the set $\{Q_1\}$ gives the set
we were looking for. Otherwise, there exists a minimal $j$ such that
$2 \le j \le t$ and $Q_1\cdot D(R_{j1})=0$. Then we
replace $\E$ by $\E \cup \{Q_1\}$ and the set
$$
\{ R_{11}, R_{12}\} \cup ...\cup \{ R_{t1}, R_{t2}\}
$$
by
$$
\{ R_{j1}, R_{j2} \mid 1\le j \le t,\  Q_1\cdot D(R_{j1})=0\}
$$
and repeat this procedure.
\enddemo

\subsubhead
5.3. Basic Theorems
\endsubsubhead

We want to prove the following basic statement.

\proclaim{Basic Theorem 2.5.8} Let $X$ be a Fano 3-fold with terminal
$\Bbb Q$-factorial singularities. Assume that $X$ does not have
a small extremal
ray, and Mori polyhedron $\nem (X)$ does not have a face of Kodaira
dimension 1 or 2.

Then we have the following statements for the $X$:

(1) The $X$ does not have a pair of extremal rays of the type $\bb_2$
(thus, in notation above, the invariant $n=0$) and Mori polyhedron
$\nem (X)$ is simplicial.

(2) The $X$ does not have more than one extremal ray
of the type (I).

(3) If $\E$ is an extremal set of $k$ extremal rays of $X$, then
the $\E$ has one of the types:
$\Aa_1\amalg (k-1)\cc_1$, $\dd_2\amalg (k-2)\cc_1$,
$\cc_2\amalg (k-2)\cc_1$,
$k\cc_1$ (we use notation of Theorem 2.3.3).

(4) We have the inequality for the Picard number of $X$:
$$
\rho (X) = \dim N_1(X) \le 7.
$$
\endproclaim

\demo{Proof} We use notations introduced in the Section 5.2.
We devide the proof on several steps.

Let us consider extremal rays
$$
\E_0 = \{ R_{11}, R_{12}\}\cup \{ R_{21}, R_{22}\}\cup ...
\cup \{ R_{n1}, R_{n2}\}.
$$
Let
$$
\E_0^{ind} =
\{ R_{11}, R_{12}\} \cup \{ R_{21}, R_{22}\}\cup ...\cup \{ R_{m1}, R_{m2}\},
$$
and
$$
\E_0^{dep} =
\{ R_{(m+1)1}, R_{(m+1)2}\} \cup \{ R_{(m+2)1}, R_{(m+2)2}\}\cup ...\cup
\{ R_{n1}, R_{n2}\}.
$$
By Lemma 2.5.4, the set $\E_0$ is extremal.
Let $\E$ be a maximal
extremal set of extremal rays which contains $\E_0$.
Let $\E_1=\E - \E_0$. By Proposition 2.3.8,
$\sharp \E_1=\rho (X)-1-\dim [\E_0]$. By Theorem 2.3.3,
for $S\in \E_1$,
the divisor $D(S)$ has no a common point with divisors
$D(R_{i1})$, $1\le i \le n$.

\proclaim{Lemma 2.5.9} Let $X$ satisfies the conditions of
Theorem 2.5.8. Let $Q$ be an extremal ray such that $Q$ is
different from extremal rays $R_{ij}$, $1\le i \le n$,
$1\le j\le 2$, and
the set $\E_1\cup \{ Q \}$ is extremal. Then the $Q$ has the type
$(II)$ and there exists exactly one $i$ such that
$1\le i \le n$ and $Q \cdot D(R_{i1})>0$ and
$D(Q) \cap D(R_{j1})=\emptyset $ if $j\not= i$.
\endproclaim

\demo{Proof} Assume that $Q$ has the type (I). Then the divisor
$D(Q)$ has no a common point with the divisors $D(R_{i1})$,
$1\le i\le n$. By Theorems 2.3.3, 2.3.6 and Lemma 2.5.1, the set
$\{Q\} \cup \E_1 \cup \E_0$ is extremal. We then get a
contradiction with the condition that $\E_1 \cup \E_0$ is a
maximal extremal set. Thus, the extremal ray $Q$ has the type (II).

If $D(Q)$ has no a common point with divisors $D(R_{i1})$, $1\le i \le n$,
we get a contradiction by the same way.
Thus, there exists $i$ such that $1\le i \le n$
and $D(Q) \cap D(R_{i1})\not= \emptyset$.
Let us consider a projectivisation
$P\nem(X)$. By Lemma 2.2.2, $P\nem (X, D(Q))$ is an interval with
two ends. Its first end is the vertex $PQ$
and its second end is a point of the edge $P(R_{i1}+R_{i2})$ of the
convex polyhedron $P\nem (X)$. Thus, the $i$ is defined by the
extremal ray $Q$. Evidently, $Q\cdot D(R_{i1})>0$.
\enddemo

\proclaim{Lemma 2.5.10} With conditions of Lemma 2.5.9 above,
assume  that $m+1 \le i \le n$. Then there exists exactly one
extremal ray $Q=Q_i$ with conditions of Lemma 2.5.9: thus,
the set $\E_1 \cup \{ Q_i\}$ is extremal and
$Q_i\cdot D(R_{i1})>0$, and $D(Q_i)\cap D(R_{j1})=\emptyset$
if $j \not= i$.
\endproclaim

\demo{Proof} The
$$
\beta =\sum_{S\in \E_1} {S} + \sum _{R\in \E_0}{R}
$$
is a face of $\nem (X)$ of highest dimension $\rho(X)-1$, and
$$
\beta_i =\sum_{S\in \E_1} {S} + \sum _{R\in \E_0-\{ R_{i1},R_{i2}\}}{R}
$$
is a face $\beta_i\subset \beta \subset \nem (X)$ of the dimension
$\rho (X)-2$ and of the codimension one in $\beta $. (Here
we use that $m+1 \le i \le m+k$). It follows that
there exists exactly one face $\beta _i^\prime$ of
$\nem (X)$ such that $\beta_i^\prime$ contains $\beta_i$,
$\dim \beta_i^\prime = \rho (X)-1$, and $\beta_i^\prime \not=\beta$.
By Theorems 2.3.3 and 2.3.6, and Lemma 2.5.9,
$\beta_i^\prime = \beta_i + Q_i$ where $Q_i$ is an extremal ray
such that the set
$\E_1\cup \{ Q_i\} \cup
(\E_0-\{ R_{i1},R_{i2}\})$ is extremal, and the
ray $Q_i$ has properties of Lemma 2.5.10. It follows that the
$Q_i$ is unique and does exist.
\enddemo

\proclaim{Lemma 2.5.11} In notation above, the set
$\E_1\cup \E_0^{ind} \cup \{ Q_{m+1},...,Q_{n} \}$
is extremal.
\endproclaim

\demo{Proof} We apply Lemma 2.5.7 to
$\E=\E_1\cup \E_0^{ind}$ and $\E_0^{dep}$. By Lemma 2.5.7,
there are extremal rays
$Q^\prime_{m+1},....Q^\prime_{m+r}$ such that the
set $\E_1 \cup \E_0^{ind}\cup
\{ Q^\prime_{m+1},....Q^\prime_{m+r}\}$
is extremal and for any $i$, $m+1\le i \le m+r$,
there exists $j$, $m+1\le j \le n$,
such that $Q^\prime_i\cdot D(R_{j1})>0$. Moreover,
for any $j$, $m+1\le j\le n$, there exists an extremal
ray $Q_i$, $m+1\le i \le m+r$, such that
$$
Q^\prime_i \cdot D(R_{j1})>0.
$$

By Lemmas 2.5.9 and 2.5.10, $r=k$ and
$\E_1 \cup \E_0^{ind}\cup
\{ Q^\prime_{m+1},....Q^\prime_{m+r}\}=
\E_1\cup \E_0^{ind} \cup \{ Q_{m+1},...,Q_{n} \}$.
\enddemo

\proclaim{Lemma 2.5.12} The set $\E_0^{dep}$ is empty.
\endproclaim

\demo{Proof} By Lemmas 2.5.9, 2.5.10 and 2.5.11, the set of
extremal rays
$U=\E_1\cup \E_0^{ind} \cup \{ Q_{m+1},...,Q_{n} \}$
is a maximal extremal set which contains
$\E_1\cup \E_0^{ind}$ and does not contain
extremal rays from $\E_0^{dep}$.
Assume that $k=n-m\not=0$. Then $k\ge 2$ and
$\dim U=\rho_0 (X) - 1 + 2m + k$. But the dimension of
a face of $\nem (X)$ of highest dimension is equal to
$\rho (X)-1 = \rho_0 (X) - 1 + 2m + k + \delta$
where $\delta \ge 1$. Thus, the extremal set $U$ is not maximal,
and there exists another extremal ray $S$ such that
$U\cup \{ S\}$ is extremal. By definition of $U$, the
$S\in \E_0^{dep}$. Let $S=R_{i1}$ where $m+1 \le i \le n$.
Since $Q_i\cdot D(R_{i1})>0$, by Theorem 2.3.3, the extremal set
$\{ Q_i, R_{i1}\}$ has the type $\cc_2$. Thus, $R_{i1}\cdot D(Q_i)=0$.
By definition of the set $\E_0^{dep}$, there exists a linear
dependence
$\sum_{l=m+1}^{l=n}{a_{l1}R_{l1}+a_{l2}R_{l2}}=0$
where $a_{i1}\not=0$ and $a_{i2}\not=0$. Multiplying the equality
above on $D(Q_i)$, we get $a_{i2}=0$. Thus, we get a contradiction.
(Compair with Lemma 2.5.6.)
\enddemo

\proclaim{Lemma 2.5.13} The set $\E_0 ^{ind}$ is empty.
\endproclaim

\demo{Proof} Since $\E_0^{dep}=\emptyset$, the set
$U=\E_1\cup \E_0^{ind}=\E_1\cup \{R_{11}, R_{12}\}
\cup...\cup \{ R_{m1}, R_{m2}\}$
is a maximal extremal set. It follows that
$U$ generates a simplicial face of $\nem (X)$ of codimension 1.
Thus,
$U_1=\E_1\cup \E_0^{ind}-\{R_{m2}\} =
\E_1\cup \{R_{11}, R_{12}\}
\cup...\cup \{ R_{(m-1)1}, R_{(m-1)2}\} \cup
\{ R_{m1}\}$ generates a simplicial face of $\nem (X)$ of codimension 2.
It follows that there exists an extremal ray $Q_{m2}$ such that
$U_1^\prime =
\E_1\cup \{R_{11}, R_{12}\}
\cup...\cup \{ R_{(m-1)1}, R_{(m-1)2}\} \cup
\{ R_{m1}\}\cup \{Q_{m2}\}$
generates a simplicial face of $\nem (X)$ of codimension 1, and
$Q_{m2}$ is different from $R_{m2}$. By Lemma 2.5.9, $Q_{m2}\cdot
D(R_{m1})>0$. Thus, by Theorem 2.3.3,
$\{Q_{m2}, R_{m1}\}$ is an extremal set of the type
$\cc_2$ where $R_{m1}\cdot D(Q_{m2})=0$.

Similarly, we can find an extremal ray $Q_{m1}$ such that
the set $\{ Q_{m1}, R_{m2}\}$ is extremal of the type $\cc_2$
where $R_{m2}\cdot D(Q_{m1})=0$. Then we get a contradiction
with Lemma 2.5.6. Thus, $m=0$, and the set $\E_0^{ind}=\emptyset$.
\enddemo

Thus, we proved that
$X$ does not have a pair of extremal rays of
the type $\bb_2$. By Theorem 2.3.3 and Proposition 2.3.8, the Mori
polyhedron $\nem (X)$ is then simplicial. Thus, we have proven the
statement (1).

\smallpagebreak

Now let us prove (2): $X$ does not have more than one extremal ray of
the type (I).
We need some definitions.

Let $P$ be a set of divisorial extremal rays. We say that $P$ is
{\it divisorially connected}
if there is no decomposition $P=P_1\cup P_2$ such
that both $P_1$ and $P_2$ are non-empty and for any $R\in P_1$ and
any $Q\in P_2$ divisors $D(R)$ and $D(Q)$ do not have a common point.
It defines {\it divisorially connected components}
of a set of extremal rays (we had used this definition before, see
definition before Theorem 2.3.3).
For example, for our case, a set of extremal
rays is extremal iff its divisorially connected components have the type
$\Aa_1, \bb_2, \cc_m$ or $\dd_2$.
Also, we can say what does it mean that two sets $P_1$ and $P_2$
of extremal rays are {\it divisorially joint}: this means that
there exist extremal rays $Q_1\in P_1$ and $Q_2\in P_2$ such that
divisors $D(Q_1)$ and $D(Q_2)$ have a common point (in particular,
this divisors or even extremlal rays $Q_1$, $Q_2$ may coincide).

We say that a set $P$ of divisorial extremal rays is
{\it single arrows connected} if for any two diffent extremal rays
$R_1, R_2\in P$ there exists an oriented path in the graph $G(P)$ with the
beginning in $R_1$ and terminal in $R_2$. This defines
{\it single arrows
connected components of $P$}. Like above, we can say what does it mean
that two sets $P_1$ and $P_2$ of divisorial extremal rays are
{\it single arrows joint}: either they have a common extremal ray, or
there exist extremal rays $Q_1, Q_1^\prime \in P_1$ and
$Q_2, Q_2^\prime \in P_2$ such that there are an oriented path joining
$Q_1$ and $Q_2$ and an oriented path joining $Q_2^\prime$ and $Q_1^\prime$
in the set $P_1\cup P_2$. For example, we can reformulate Lemma 1.1:
for our situation,
an $E$-set $\La$ is single arrows connected, two $E$-sets $\La$ and $\M$ are
single arrows joint (using Lemma 2.5.1, Theorem 2.4.1 and the statement (1) we
have proven).

By Lemma 2.2.2, divisors of different extremal rays of
the type (I) do not have a common point.
By Theorem 2.3.6, any set of extremal rays of
the type (I) generates a simplicial face of $\nem (X)$ of Kodaira
dimension 3. It follows that the set of extremal rays of the type (I) is
finite.
Let
$$
\{ R_1,...,R_s\}
$$
be the hole set of divisorial extremal rays of the type (I) on $X$.
We should prove that $s\le 1$.

Let $\E$ be a maximal extremal set of extremal rays on $X$
containing $\{ R_1 ,..., R_s\}$.
Let $T_1,...,T_s$ are divisorially
connected components
of $\E$ containing $R_1,...,R_s$ respectively. By Theorem 2.3.3,
these connected components are different. The connected component
$T_i$ either contains one ray $R_i$ (has the type $\Aa_1$) or
contains two rays: the ray $R_i$ and another ray of the type (II). Then
this connected component has the type $\dd_2$. Let
$$
\E_0=\E - (T_1\cup ... \cup T_s).
$$
Thus, extremal sets $\E_0,\ T_1,\ ...,\ T_s$ are divisorially disjoint.

By \cite{Ka1} and \cite{Sh}, any face of $\nem (X)$ is
contractible, and by our condition it has Kodaira dimension 3.
By Proposition 2.2.6, for any $1 \le i \le s$, there exists
an effective divisor $D(T_i)$ which is a linear
combination of divisors of rays from $T_i$ with positive coefficients,
such that $R\cdot D(T_i)<0$ for any $R\in T_i$.

Using the divisors $D(T_i)$, similarly to Lemma 2.5.7,
we can find extremal rays
$$
\{Q_1,...,Q_r\}
$$
with properties:

(a) $r \le s$;

(b) For any $i$, $1\le i \le r$, there exists $j$, $1\le j \le t$,
such that $Q_i\cdot D(T_j)>0$ (in particular, $Q_i$ is different
from extremal rays of $\E$ and do not have the type (I));

(c) For any $j$, $1\le j\le s$, there exists an extremal
ray $Q_i$, $1\le i \le r$, such that
$$
Q_i \cdot D(T_j)>0;
$$

(d) The set
$\E_0 \cup \{ Q_1,...,Q_r\}$
of extremal rays
is extremal, and extremal rays
$\{ Q_1,...,Q_r\}$ are linearly independent.

By our condition, all extremal rays on $X$ are divisorial. Thus,
by (b), the extremal rays $Q_1,...,Q_r$ have the type (II).

Let us take a ray $Q_i$ and let $Q_i\cdot D(T_j)>0$.
By Theorem 2.3.3, the
$T_j$ generates a simplicial face $\gamma_j$
of $\nem (X)$. Since $T_j$ has the
type $\Aa_1$ or $\dd_2$, one can see easily using Lemma 2.2.3, that
any curve of divisors of rays from $T_j$ belongs to this face.
It follows that $\nem (X, D(Q_i))$ is a 2-dimensional angle
bounded by the ray $Q_i$ and a ray from the face $\gamma_j$
since the divisor $D(Q_i)$ evidently has a common curve with
one of divisors $D(R)$ of $R\in T_j$. Since any two sets from
$T_1,...,T_s$
do not have a common extremal ray, the faces
$\gamma_1,...,\gamma_s$ do not have a common ray (not necessarily
extremal).
It follows that the angle $\nem (X, D(Q_i))$ does not have a common ray
with the face $\gamma_k$ for $k\not=j$.
Thus, the divisor $D(Q_i)$ do not have a
common point with divisors of rays $T_k$. It follows that
$r=s$ and we can choose
numeration $Q_1,...Q_s$ such that $Q_i\cdot D(T_i)>0$ but
$D(Q_i)$ do not have a common point with divisors of extremal rays
$T_j$ if $i\not=j$ (the $Q_i$ is divisorially disjoint with
$T_j$ if $i\not=j$).

Let us fix $i$, $1\le i \le n$. The set
$\E\cup \{Q_i\}$ is not extremal since $\E$ is maximal extremal.
Since $Q_i$ is divisorially disjoint with divisorially connected
components  $T_j$ of $\E$ if $j\not=i$, it follows that
the set $\E_0 \cup T_i \cup \{ Q_i\}$
is not extremal also.
Then, the last set contains  an $E$-set $\La_i$.
Since $\E_0\cup T_i$ is  extremal,
the set $\La_i$ contains the ray $Q_i$. The
$\La_i$ is single arrows connected (we have mentioned this above).
Let $\E_{0i} = \La_i \cap \E_0$. We claim that the subset
$$
U_i=\E_{0i}\cup \{ Q_i\}
$$
of $\La_i$ is single arrows connected too. Really, the  set
$\La_i^\prime =\La_i-(\E_{0i}\cup \{ Q_i\} )$
is not empty since the set
$\E_{0i}\cup \{ Q_i\}$ is extremal as a subset of the extremal
set $\E_0\cup \{Q_1,...,Q_s\}$.
Since $\La_i$ is single arrows connected, there exists a
shortest oriented path
in $\La_i$  connecting any $S\in \E_{0i}$ with $\La_i^\prime$.
But $\E_{0i}$ and $\La_i^\prime$ are divisorially (and then
single arrows) disjoint. Thus, the beginning of this
path is a path in $\E_{0i}\cup \{ Q_i\} $ joining
$S$ with the ray $Q_i$. The same considerations prove that
there exists an oriented path  in $\E_{0i}\cup \{ Q_i\}$
joining $Q_i$ with $S$.

By Lemma 1.1, for $i\not= j$, the sets $\La_i$ and $\La_j$ are single
arrows joint. Using the fact that
$\La_i^\prime$ and $\La_j^\prime$ are divisorially disjoint, one
similarly can prove that the sets
$U_i = \E_{0i}\cup \{ Q_i\} $ and
$U_j = \E_{0j}\cup \{ Q_j\} $ are single arrows joint.

Let $U$ be the union of all sets
$U_i = \E_{0i}\cup \{ Q_i\} $, $i=1,..., s$.
Then $U$ is single arrows connected
extremal set of extremal rays of
the type (II).
By Theorem 2.3.3, it is possible only if $U$ is either empty or
has the type $\cc_1$. It follows that $s\le 1$.
This proves the statement (2).

\vskip5pt

Let us prove (3).

We use the following

\proclaim{Statement} The contraction
of a ray $R$ of the type (II) on $X$
gives a Fano 3-fold $X^\prime$ with terminal
$\Bbb Q$-factorial singularities and
without small extremal rays and without
faces of Kodaira dimension 1 or 2 for $\nem (X^\prime )$.
Extremal sets $\E^\prime$ on $X^\prime$ are in one to one
correspondence with extremal sets $\E$ on $X$ which contain the ray $R$.
\endproclaim

\demo{Proof} Let $\sigma :X\to X^\prime $
be a contraction of $R$.
The $X^\prime$ has terminal $\Bbb Q$-factorial
singularities by \cite{Ka1} and \cite{Sh}.
We have,
$K_X = \sigma ^*(K_{X^\prime}) + dD(R)$. Multiplying this
equality on $R$ and using Proposition 2.3.2, we get that
$d=1$. By the statement (1), it follows that
$\sigma ^*( - K_{X^\prime}) = -K_X+D(R)$ is nef and only contracts the
extremal ray $R$. Then $-K_{X^\prime}$ is ample on $X^\prime$ and
$X^\prime$ is a Fano 3-fold with terminal $\Bbb Q$-factorial singularities.
Faces of $\nem (X^\prime )$ are in one to one correspondence with faces of
$\nem (X)$ which contain the $R$. Contractions of faces of $\nem (X^\prime)$
are dominated by contractions of the corresponding faces of $\nem (X)$.
It follows the last statement.
\enddemo

Let $\E=\{ R_1,...,R_k\} $ be an extremal set on $X$. By Theorem 2.3.3,
it has connected components of the type $\Aa_1, \bb_2, \cc_m$ or $\dd_2$.
Moreover, by (1) and (2), it does not have a connected component of the type
$\bb_2$ and does not have more than one connected component of the type
$\Aa_1$. By Statement above, the same should be true for the extremal set
$\E^\prime$ one gets by the contraction of any extremal ray $R_i$ of the
type (II) of $\E$. It follows the statement (3).

Now we prove (4): $\rho (X)\le 7$.

First, we show how to prove $\rho (X)\le 8$ applying Theorem 1.2 to the face
$\gamma=\M(X)$ of $\dim \M(X) =m=\rho (X)-1$. By the statement (1) of
Theorem 2.5.8 and Theorems 2.3.3 and 2.4.1, the $\M(X)$ is simple and
all conditions of Theorem 1.2 are valid for some constants $d, C_1, C_2.$
By Theorem 2.4.1, we can take $d=2$.
By the proof of Theorem 1.2, we should find the constants $C_1$ and $C_2$
for maximal extremal sets $\E$ only (only this sets we really use). Thus,
$\sharp \E =m$. By the statement (3), then the constants
$C_1\le 2/m$ and $C_2 = 0$. Thus, we get
$m<(16/3)2/m+6$. Then, $m=\rho (X)-1\le 7$, and $\rho(X)\le 8$.

To prove the better inequality $\rho (X)\le 7$, we should analyze the proof of
Theorem 1.2 for our case more carefully.
We will show that the conditions
of Lemma 1.4 hold for the $\M(X)$
with the constants $C=0$ and $D=2/3$.
By Lemma 1.4, we then get the inequality $\rho (X)\le 7$ we want to prove.

Like for the proof of Theorem 1.2, we introduce a weight of an
oriented angle, but using a new
formula:
$\sigma (\angle)=2/3$ if $\rho (R_1(\angle),R_2(\angle))=1$, and
$\sigma (\angle)=0$ otherwise.

By the statement (3) of Theorem 2.5.8, the condition (1) of
Lemma 1.4 holds with constants $C=0$ and $D=2/3$.

Let us prove the condition (2) of Lemma 1.4. For $k=3$ (triangle) it is
true since an $E$-set which has at least 3 elements has the type (a) of
Theorem 2.4.1 (see the proof of Theorem 1.2). Thus,
the triangle has at least three oriented angles with the
weight $2/3$.  For $k=4$ (quadrangle),
we proved (when we were proving Theorem 1.2) that one can find at least two
oriented angles of the quadrangle such that
for any of them $\rho (R_1(\angle), R_2(\angle ))$ is finite. By the
statement (3) of Theorem 2.5.8, then
$\rho (R_1(\angle), R_2(\angle ))=1$. Thus, the quadrangle has at least
two oriented angles of the weight $2/3$.
This finishes the proof of Theorem 2.5.8.

\enddemo

Now, we give an application
of the statement (2) of Theorem 2.5.8 to geometry of Fano 3-folds.

Let us consider a Fano 3-fold $X$ and blow-ups $X_p$ in different
non-singular points $\{x_1,...,x_p\} $ of $X$. We say that this is
a Fano blow-up if $X_p$ is Fano.  We have the following very simple

\proclaim{Proposition 2.5.14}
Let $X$ be a Fano 3-fold with terminal $\Bbb Q$-factorial
singularities and without small extremal rays.
Let $X_p$ be a Fano blow up of $X$.

Then for any small extremal ray $S$ on $X_p$, the
$S$ has a non-empty intersection with one of exceptional
divisors $E_1,...,E_p$ of this blow up and does not belong to any of
them. Moreover, the exceptional divisors $E_1,...,E_p$ define $p$
extremal rays $Q_1,...,Q_p$ of the type (I) on $X_p$ such that
$E_i=D(Q_i)$.
\endproclaim

\demo{Proof} The last statement is clear.
Let $S$ be a small extremal ray on $X_p$ which does not intersect
divisors $E_1,...,E_p$. Let $H$ be a general $nef$ element orthogonal
to $S$. Let $l_1,...,l_n$ are lines  which generate extremal rays
$Q_1,...,Q_p$. Then the divisor
$H^\prime = H+(l_1\cdot H)/(-l_1\cdot E_1)E_1+...
+(l_p\cdot H)/(-l_p \cdot E_p)E_p$
is a nef divisor on $X_p$
orthogonal to all extremal rays $Q_1,...,Q_p, S$, and
$(H^\prime)^3>H^3>0$.  This proves that the extremal rays
$Q_1,...,Q_p, S$ generate a face of $\nem (X_p)$ of Kodaira
dimension $3$. Then, by contraction of extremal rays
$Q_1,...,Q_p$, the image of $S$ gives a small extremal ray on $X$.

This gives a contradiction.
\enddemo

It is known that a contraction of a face of
Kodaira dimension 1 or 2 of $\nem (Y)$ of a Fano 3-fold $Y$ has
a general fiber which is a rational surface or curve respectively,
because this contraction has relatively negative canonical class.
See \cite{Ka1}, \cite{Sh}.
It is known that a small extremal ray is rational \cite{Mo2}.

Then, using the statement (2) of
Theorem 2.5.8 and Proposition 2.5.14, we can divide
Fano 3-folds of Theorem 2.5.8 on the following 3 classes:

\proclaim{Corollary 2.5.15}
Let $X$ be a Fano 3-fold with terminal $\Bbb Q$-factorial
singularities and without small
extremal rays, and without faces of Kodaira dimension 1 or 2 for
Mori polyhedron. Let $\epsilon$ be the number of extremal rays
of the type (I) on $X$ (by Theorem 2.5.8, the $\epsilon \le 1$).

Then there exists $p$, $1\le p \le 2-\epsilon$, such that
$X$ belongs to one of classes (A), (B) or (C) below:

(A) There exists a Fano blow-up $X_p$ of $X$ with a face of Kodaira
dimension 1 or 2. Thus, birationally, $X$ is a fibration on
rational surfaces over a curve or rational curves over a surface.

(B) There exist Fano blow-ups $X_p$ of $X$ for
general $p$ points on $X$ such that
for all these blow-ups the
$X_p$ has a small extremal ray $S$. Then images of curves of $S$ on
$X$ give a system of rational curves on $X$ which cover a Zariski
open subset of $X$.

(C) There do not exist Fano blow-ups $X_p$ of $X$ for general
$p$ points.

We remark that for Fano 3-folds with Picard number $1$ the
 $\epsilon =0$. Thus, $1\le p \le 2$.
\endproclaim

We mention that the statements (3) and (4) of Theorem 2.5.8 give similar
information for blow ups of $X$ at curves. Of course, it is more
difficult to formulate these statements.

\enddocument

\end